\documentclass[smallextended]{svjour3}  
\usepackage{latexsym,amsmath,amssymb,theorem}
\usepackage[all]{xy}

\usepackage[usenames,dvipsnames]{pstricks}
\usepackage{epsfig}

\usepackage{color} 
\usepackage{hyperref}
\usepackage{mathbbol}
\usepackage{marginnote} 
\theoremstyle{change} 

\newenvironment{Proof}{\noindent \bf Proof:   \rm}{\hspace*{\fill}
$\square$ \medskip}

{\theorembodyfont{\itshape}
\newtheorem{Lemma}{Lemma}
\newtheorem{Proposition}[Lemma]{Proposition}
\newtheorem{Fact}[Lemma]{Fact}

\newtheorem{Corollary}[Lemma]{Corollary}
\newtheorem{Theorem}[Lemma]{Theorem}}
{\theorembodyfont{\upshape}

\newtheorem{Definition}[Lemma]{Definition}
\newtheorem{Remark}[Lemma]{Remark}
\newtheorem{Remarks}[Lemma]{Remarks}

\newtheorem{Examples}[Lemma]{Examples}}

\newcommand{\tr }{\textcolor{red}}

\newcommand{\ot}{\otimes}

\newcommand{\N}{\mathbb{N}}

\newcommand{\BC}{\mathbf{C}}

\newcommand{\BD}{\mathbf{D}}

\newcommand{\BB}{\mathbf{B}}
\newcommand{\BA}{\mathbf{A}}

\newcommand{\BS}{\mathbf{S}}

\newcommand{\ring}{\mathbf{Ring}}
\newcommand{\Coring}{\mathbf{Coring}}
\newcommand{\Mod}{\mathbf{Mod}}

\newcommand{\Comon}{\mathbf{Comon}}
\newcommand{\Bimon}{\mathbf{Bimon}}

\newcommand{\op}{\mathsf{op}}

\newcommand{\id}{\mathsf{id}}
\newcommand{\Id}{\mathsf{Id}}

\newcommand{\ma}{\mathsf{A}}
\newcommand{\mb}{\mathsf{B}}
\newcommand{\mc}{\mathsf{C}}

\newcommand{\mh}{\mathsf{H}}
\newcommand{\mk}{\mathsf{K}}
\newcommand{\sm}{\mathsf{M}}
\newcommand{\si}{\mathsf{I}}
\newcommand{\sn}{\mathsf{N}}
\newcommand{\sfp}{\mathsf{P}}
\newcommand{\Set}{\mathbf{Set}}

\newcommand{\Alg}{\mathbf{Alg}}

\newcommand{\Mon}{\mathbf{Mon}}

\newcommand{\Vect }{\mathbf{Vect}_k}
\newcommand{\Bialg}{\mathbf{Bialg}}
\newcommand{\Coalg}{\mathbf{Coalg}}
\newcommand{\Hopf}{\mathbf{Hopf}}

\newcommand{\lra}{\rightarrow}

\newcommand{\xra}{\xrightarrow}

\begin{document}

  \title{Generalizations of the Sweedler dual}
  
 \author{ 
 Hans-E. Porst 
 \and  
 Ross Street\thanks{The second author gratefully acknowledges the support of Australian Research Council Discovery Grant DP130101969.}
}
\institute{Hans-E. Porst \at
Department of Mathematical Sciences, University of Stellenbosch, \\  Stellenbosch, South Africa \\
\email{porst@math.uni-bremen.de} \\
\emph{Permanent address: Department of Mathematics, University of Bremen, 28359 Bremen, Germany}
\and
 Ross Street \at
Centre of Australian Category Theory, Department of Mathematics, \\  Macquarie University NSW 2109, Australia \\
\email{ross.street@mq.edu.au}
}
\date{}
\dedication{To the memory of Horst Herrlich}

\maketitle
\begin{abstract}
As  left adjoint to the dual algebra functor, Sweedler's finite dual construction is an important tool in the theory of Hopf algebras over a field. 
We show in this note that the left adjoint to the dual algebra functor, which exists over arbitrary  rings,  shares a number of properties with the finite dual. Nonetheless the requirement that it should map Hopf algebras to Hopf algebras needs the extra assumption that this left adjoint should map an algebra into its linear dual.  
We identify a condition guaranteeing that Sweedler's construction works when generalized to noetherian commutative rings. We establish the following two apparently previously unnoticed dual adjunctions: 
For every commutative ring $R$ the left adjoint of the dual algebra functor on the category of $R$-bialgebras has a right adjoint. 
This dual adjunction can be restricted to a dual adjunction on the category of Hopf $R$-algebras,
provided that $R$ is noetherian and absolutely flat.


\end{abstract}


\section*{Introduction}
Given an $R$-coalgebra $\mc$ over a commutative unital ring $R$, there is a \em dual algebra \em $\mc^\ast$ whose underlying $R$-module is the linear dual $C^\ast$ of the underlying module  $C$ of $\mc$. 
Conversely, if the underlying module $A$ of the $R$-algebra $\ma$ is finitely generated projective, there is a \em dual coalgebra \em $\ma^\ast$ on the linear dual of $A$. 
In this case the linear isomorphism $A\simeq ({A^\ast})^\ast$ induces an isomorphism $\ma\simeq ({\ma^\ast})^\ast$.  These isomorphisms then constitute a duality between the categories	${_{f\!gp}\Coalg_R}$ and ${_{f\!gp}\Alg_R}$ of finitely  generated projective coalgebras and algebras respectively. These facts are well known.

In the case where $R=k$ is a field, Sweedler's \em finite dual coalgebra construction $\ma^\bullet$ \em \cite{Sweedler} of a $k$-algebra $\ma$ provides an extension of the above construction to arbitrary algebras.
The underlying vector space $\ma^\circ$ of $\ma^\bullet$ is the subspace of $A^\ast$ consisting of those linear forms whose kernel contains a cofinite ideal.  
This extension is no longer an inverse equivalence to the (restricted) dual algebra	functor, but only a left adjoint. 
In addition, $\ma^\bullet$ underlies a Hopf algebra if $\ma$ does.

Attempts have been made to generalize this to more general commutative rings than just fields. 
In \cite{Wis} and \cite{Cao} the authors have shown that  Sweedler's original construction of $\ma^\circ$, which clearly makes sense over any commutative ring $R$, can be made into an $R$-coalgebra by the same (see \cite{Cao}) or a related (see \cite{Wis}) comultiplication as in the field case for larger classes of rings. 
They then succeeded in showing that their constructions also send $R$-Hopf algebras to $R$-Hopf algebras. 
The questions whether  their constructions are  functorial or even the construction of a left adjoint 
to the dual algebra functor are not considered there.

On a different note it has been shown in \cite{HEP_AJSE}  that the dual algebra functor indeed has a left adjoint, for every commutative ring $R$. This clearly is Sweedler's dual coalgebra functor, if $R=k$ is a field. Nothing is known however, to what extent this left adjoint shares any of the further properties of Sweedler's construction, when $R$ is not a field.

It is the purpose of this paper to close the missing gaps. 
We will prove that the constructions of \cite{Wis} and \cite{Cao} indeed provide left adjoints.  
We will show as well, using methods from the theory of monoidal categories, that, for every commutative ring $R$, the left adjoint $\mathsf{Sw}$ of the dual algebra functor $D = (-)^\ast$, which we call the \em generalized Sweedler dual  functor\em, shares certain properties with Sweedler's construction but fails to map Hopf algebras to  Hopf algebras in general. 

In more detail: 
\begin{enumerate}
\item  Instead of an embedding of $\ma^\circ$ into $A^\ast$ there is only a  canonical map $\kappa_\ma$ from the underlying algebra of $\mathsf{Sw}(\ma)$ into the linear dual of $A$. 
\item The generalized Sweedler dual  functor $\mathsf{Sw} $, as with the ordinary finite dual, extends the dual coalgebra functor from ${_{f\!gp}\Alg_R}$ to ${\Alg_R}$.
\item  The generalized Sweedler dual functor is merely opmonoidal, whereas the ordinary finite dual functor is \em strong \em (op)monoidal.  
\item $\mathsf{Sw}(\ma)$ is a bialgebra when $\ma$ is the underlying algebra of a bialgebra. 
Yet $\mathsf{Sw} $ maps only those Hopf algebras $\mh$ to Hopf algebras in the case where the canonical map $\kappa_{\mh^a}$ of the underlying algebra of $\mh$ is an embedding.
\end{enumerate}

This leads to the question of how to construct a generalized finite dual such that $\kappa$ is a monomorphism. 
Categorically analyzing Sweedlers's original construction defined for arbitrary rings literally as in the field case,
we can not only motivate this construction, but also identify its obstacles for arbitrary rings. We develop a criterion for this construction to support a functor and eventually the desired left adjoint. By showing that the constructions of \cite{Wis} and \cite{Cao} meet this criterion's hypotheses, we get the result mentioned earlier, that their constructions indeed provide left adjoints to the dual algebra functor.

As a by-product, the categorical methods used enable us to establish two hitherto unknown dual adjunctions. 
For every commutative ring $R$, the generalized finite dual
$\mathsf{Sw} $, considered as a functor $\Bialg_R\rightarrow\Bialg_R^\op$ has a right adjoint. 
 This dual adjunction can be restricted to  a dual adjunction on $\Hopf_R$, provided that $R$ is  
 noetherian and absolutely flat.  

Concerning the presentation we would like to add the following remarks: Some of our results are of interest  for not  necessarily commutative rings as for example the fact that Sweedler's dual $R$-ring functors have left adjoints, and these extend the familiar duality for finitely generated projective $R$-rings and corings, respectively (see e.g. \cite{Swe} or \cite{LP_HP}).  We, thus, assume rings $R$ to be commutative only where this is needed, that is, from subsection \ref{subsec:com} onwards.  Moreover, since the dual algebra functor is a special instance of a functor $\bar{G}\colon\Mon\BC\lra\Mon\BD$  between categories of monoids in monoidal categories, induced by a monoidal functor $G\colon\BC\lra\BD$, and since some of our results are of interest even in this general setting, as for example the fact that such a functor $\bar{G}$ is symmetric monoidal, provided that $G$ is, we start our investigations there (that is, this level of abstraction is not chosen merely for the sake of utmost generality).

 The paper is organized as follows: 
 
Section \ref{sec:basics} recalls the necessary tools concerning monoidal  functors, in particular the  concepts  of \em semi-dualization functors \em and \em dual monoid functors \em in monoidal closed categories.

 Section 2 is devoted to properties of the left adjoint $F$ of a functor induced by a monoidal functor $G\colon\Mon\BC\lra\Mon\BD$. Starting with formulating necessary conditions on the algebraic structure of the $\BC$-monoid $F\ma$ for some $\BD$-monoid $\ma$, subsection \ref{ssec:struct} ends with a criterion for the existence of a functor $F\colon\Mon\BD\lra\Mon\BC$ satisfying these. Subsection \ref{ssec:ext} is concerned with the extension property 2 above, while in subsection \ref{subsec:com} it will be shown that each such $F$ is opmonoidal and thus maps bimonoids to bimonoids if $G$ is symmetric monoidal, thereby proving statements 3 and 4.

 Section \ref{sec:con} deals with generalizations of the Sweedler construction to more general rings.
We first sharpen the criterion formulated in Section \ref{ssec:struct} to an adjoint functor theorem, which can be applied in particular to certain dual comonoid functors over noetherian rings. In subsection \ref{ssec:dualcomonoid}
we discuss a method generalizing Sweedler's construction $(-)^\circ$ subject to a choice of a subcategory $\BS$ of $\Mod_R$ and explore, under what conditions this can be lifted to a generalized finite dual  functor by the criterion just mentioned. 
These results not only allow for strenghtened versions of the dual coalgebra constructions of \cite{Wis} and \cite{Cao}, as finally shown in Section \ref{sec:fin}, but also
 contribute to a better understanding of this construction. 
 We conclude with some comments on the choice of $\BS$ and on generalizations.

\section{Preliminaries}\label{sec:basics}
\subsection{Monoidal categories and functors}\label{ssec:mon_cats}
Throughout $\BC = (\BC, -\otimes -, I)$ denotes a  monoidal category. 
If $\BC$ is moreover symmetric monoidal, the symmetry will be denoted by $\sigma = (C\otimes D\xrightarrow{\sigma_{CD}} D\otimes C)_{C,D}$. 

By $\BC^\op$ we denote the dual of $\BC$ with tensor product and unit as in $\BC$, while $\BC^t$ denotes the transpose of $\BC$ with tensor product $C\ot^tD = D\ot C$.

Recall \cite{StBk}\footnote{We make this choice so that an object with a right dual will then have a right internal hom. The opposite choice appears in \cite{LP_HP}.}that  $\BC$ is called \em monoidal left closed \em provided that, for each $\BC$-object $C$ the functor $C\otimes -$ has a right adjoint $[C,-]_l$. 
If each functor  $-\otimes C$ has a right adjoint, denoted by $[C,-]_r$, $\BC$ is called \em monoidal right closed\em.  If $\BC$ is monoidal left and right closed, it is called \em monoidal closed\em. 
The counits $C\otimes [C,X]_l\rightarrow X$ and $ [C,X]_r\otimes C\rightarrow X$ of these adjunctions will be denoted by $ev^l$ and $ev^r$ respectively. 

By parametrized adjunctions (see \cite{MacL}) one thus has functors $[-,-]_r$ and $[-,-]_l$ \ $\BC^\op\times \BC\rightarrow\BC$. 
In particular, for each	$X$ in $\BC$, there are the contravariant functors $[-,X]_r$ and $[-,X]_l$ on $\BC$.
For $C\xrightarrow{f}D$ in $\BC$, the morphism $[D,X]_l\xrightarrow{[f,X]_l}[C,X]_l$ is the unique morphism such that the following diagram commutes. 
$$\begin{minipage}{6cm}
\xymatrix@=2.5em{
 C\otimes[D,X]_l    \ar[r]^{ C\otimes [f,X]_l }\ar[d]_{ f\otimes[D,X]_l }&{C\otimes [C,X]_l } \ar[d]^{ev_l }  \\
   D\otimes[D,X]_l    \ar[r]_{ ev_l}          &    X}
\end{minipage}$$
Similarly for $[-,X]_r$. 
If $\BC$ is  monoidal closed, then so is $\BC^t$; its 
internal hom-functors can be chosen as $[C,-]^t_l = [C,-]_r$ and   $[C,-]^t_r = [C,-]_l $.

We denote the categories of monoids $\sm = (M,M\ot M\xrightarrow{
m}M, I\xrightarrow{e}M)$  and of comonoids $\mc =(C, C\xrightarrow{\Delta}C\ot C,C\xrightarrow{\epsilon}I)$ in $\BC$
by $\Mon \BC$ and $\Comon\BC$, respectively.  Obviously 
$\Mon(\BC^t) = \Mon\BC \text{\  and \ } \Mon\BC^\op = (\Comon\BC)^\op$.

If $\BC =(\BC,-\ot-, I) $ is a monoidal category and $\BA$ is a subcategory of $\BC$ not necessarily closed under the monoidal structure, we denote by $\Mon\BA$ the full subcategory of $\Mon\BC$ spanned by all $\BC$-monoids $(M,m,e)$ with $M\in\BA$.

A  monoidal functor from $\BC$ to $\BC'$ will be denoted by $G = (G,\Gamma,\gamma)$ with endofunctor  
$G\colon\BC\rightarrow\BC'$, multiplication $\Gamma$  which is a natural transformation $G_1\Rightarrow G_2$ between the functors 
$G_1 = \BC\times\BC\xra{G\times G}\BD\times\BD\xra{-\otimes -}\BD \text{\ \ and \ \ } G_2 = \BC\times\BC\xra{-\otimes -}\BC\xra{G}\BD$,
and unit morphism $\gamma\colon I' \rightarrow  GI$. 

A  monoidal functor is called \em strong monoidal \em if $\Gamma$ and $\gamma$ are isomorphisms, \em strict monoidal \em if $\Gamma$ and $\gamma$ are identities, and \em normal \em if $\gamma$ is an isomorphism. 

An \em opmonoidal functor \em from  $\BC$ to $\BC'$ is a monoidal functor from $\BC^\op$ to $\BC'^\op$.\\

Given monoidal functors $F=(F,\Phi,\phi)$ and $G=(G,\Gamma,\gamma)$ from $\BC$ to $\BD$, a natural transformation $\mu\colon F\Rightarrow G$  is called a \em monoidal transformation\em, if the following diagrams commute for all $C,D$ in $\BC$.
\begin{equation*}
\begin{minipage}{6cm}
\xymatrix@=3em{
   FC\otimes FD\ar[r]^{\mu_C\otimes\mu_D}  \ar[d]_{ \Phi_{C,D} }&{ GC\otimes GD} \ar[d]^{ \Gamma_{C,D}}  \\
  F(C\otimes D)  \ar[r]_{ \mu_{C\otimes D}}          &    G(C\otimes D) 
}
\end{minipage}\hspace{2cm}\begin{minipage}{6cm}
\xymatrix@=3em{
  I \ar[dr]_{\gamma}\ar[r]^{\phi}&{ FI} \ar[d]^{ \mu_I} \\
          &    GI 
}
\end{minipage}
\end{equation*}

\begin{Remark}\label{rem:day}\rm
The following is essentially due to B. Day \cite{Day} and justifies the above terminology	for the components of a monoidal functor.

Let $\BC$ and $\BD$ be symmetric monoidal categories where $\BD$ is monoidally cocomplete; that is, {$\BD$ is cocomplete and all functors $D\ot -$  preserve colimits}. 
Then
\begin{enumerate}
\item The functor category $[\BC,\BD]$ is symmetric monoidal when equipped with the  \em Day convolution product \em $-\star -$.
\item Morphisms $\Lambda\colon F\star G \rightarrow H$ in this category  are in bijection with  families of morphisms $FC_1\otimes GC_2\rightarrow H(C_1\ot C_2)$, natural in $C_1$ and  $C_2$.
\item  Given a triple $G:=(G,\Gamma,\gamma)$, where $G$ is a  functor $\BC\rightarrow\BD$, $\Gamma$ is a natural	transformation $G_1\Rightarrow G_2$ (thus, a morphism $G\star G\rightarrow G$ in $[\BC,\BD]$), and $\gamma\colon I\rightarrow GI$ is a $\BD$-morphism, then  $G$ is a  monoidal functor $\BC\rightarrow\BD$ if and only if $G$ is a  monoid in $[\BC,\BD]$.
\item  For such triples 
$F$ and $G$ a natural transformation $F\Rightarrow G$ is a monoidal transformation if and only if it is a  monoid morphism.
\end{enumerate}
In other words: up to an equivalence of categories, monoidal functors $(G,\Gamma,\gamma)\colon \BC\xrightarrow{}\BD$ `are'  monoids $(G,\Gamma,\gamma)$ in $([\BC,\BD], -\star- )$ and monoidal transformations `are' monoid morphisms. 
This interpretation is  in fact valid for any $\BD$ in view of the fact that it can be embedded into a monoidally cocomplete category (see \cite{Day} as well).
\end{Remark}

\begin{Remark}[\cite{Kelly},\cite{KS},\cite{Ag}]\label{rem:monfcts} 
\begin{enumerate}
\item The composite of monoidal functors is a monoi\-dal functor. 
\item Monoidal functors map monoids to monoids.  {This is a consequence	of 1., since a monoid in $\BC$ is a monoidal functor $\mathbf{1}\rightarrow \BC$ with $\mathbf{1}$ the terminal category.}

In more detail: 
Let  $G:=(G,\Gamma,\gamma)\colon \BC\rightarrow\BD$ be a  monoidal functor. Then, for any monoid $(M,m,e)$ in $\BC$,
$$
 (GM, GM\otimes GM\xrightarrow{\Gamma_{M,M}}G(M\otimes M)\xrightarrow{Gm}GM, I\xrightarrow{\psi}GI\xrightarrow{Ge}GM)$$ is a monoid in $\BD$ and this construction, acting as $G$ on monoid morphisms, is functorial. 
 We call the resulting functor 
$\Mon\BC\rightarrow\Mon\BD$ the \em  functor induced by $G$ \em and denote it, with slight abuse of notation, by $G$ as well.

\item Given monoidal functors $F,G\colon\BC\rightarrow\BD$ and a monoidal transformation
 $\Gamma\colon F\Rightarrow G$, there is a natural transformation, also denoted by $\Gamma$, between the induced functors given by ${\Gamma}_{\sm,\sn} = \Gamma_{M,N}$ for monoids $\sm$ and $\sn$ in $\Mon\BC$.
 \end{enumerate}
\end{Remark}
\begin{Remark}\label{rem:opmon}
\begin{enumerate}
\item If  $(G,\Gamma,\gamma)$ is a monoidal functor and $F$ is left adjoint to $G$ (with unit $\eta$ and counit $\epsilon$), then $(F,\Phi,\phi)$ is opmonoidal, where $\phi$ corresponds under adjunction to $\gamma$ and $\Phi_{C,D}$ corresponds under adjunction to $\Gamma_{FC,FD}\circ (\eta_C\otimes\eta_D)$. We call $(F,\Phi,\phi)$ the \em opmonoidal left adjoint \em  of $(G,\Gamma,\gamma)$. The connection between a monoidal functor and its opmonoidal left adjoint is thus characterized by commutativity of the diagrams 
\begin{equation}\label{eqn:mate1}
\begin{aligned}
\xymatrix@=3em{
  C\ot D   \ar[r]^{\eta_{C\ot D}  }\ar[d]_{ \eta_C\ot\eta_D }&{GF(C\ot D) } \ar[d]^{G\Phi_{C,D} }  \\
 GFC\ot GFD   \ar[r]_{ \Gamma_{FC,FD}}          &   G(FC\ot FD) 
}\qquad
\xymatrix@=3em{
   FI  \ar[dr]^{\phi  }\ar[d]_{ F\gamma } & \\
  FGI\ar[r]^{\epsilon_I}  &    I
 }
 \end{aligned}
\end{equation}
or, equivalently	
\begin{equation}\label{eqn:mate2}
\begin{aligned}
\xymatrix@=3em{
F(GC\ot GD)   \ar[r]^{\Phi_{GC, GD}  }\ar[d]_{ F\Gamma_{C,D} }&{FGC\ot FGD } \ar[d]^{\epsilon_C\ot\epsilon_D }  \\
FG(C\ot D)   \ar[r]_{ \epsilon_{C\ot D}}          &   C\ot D
}\qquad
\xymatrix@=3em{
   I  \ar[r]^{\eta_I  }\ar[dr]_{ \gamma }&{GFI } \ar[d]^{G\phi }  \\
    &    GI
 }
 \end{aligned}
\end{equation}
\item  If  $F$ is a left adjoint of a monoidal functor and its opmonoidal structure is strong, then the unit and the counit of the adjunction are monoidal transformations;  such adjunctions are called  \em monoidal adjunctions\em. The left adjoint of a monoidal adjunction is always strong.
\end{enumerate}
\end{Remark}

\subsection{Semi-dualization in monoidal closed categories}\label{sec:dualmoncat}

Let $\BC$ be a monoidal closed category.
We call the functor $[-,I]_l$, introduced in Section~\ref{ssec:mon_cats} above, the \em left semi-dualization functor of $\BC$\em. Analogously there is the \em right semi-dualization functor \em $[-,I]_r$.

\begin{Proposition}[\cite{LP_HP}]\label{prop:dual-mon}
Let $\BC$ be a  monoidal closed category. Let  $\BC\xra{[-,I]_r}\BC^\op$ be the right and  $\BC^\op\xra{[-,I]_l}\BC$ the left semi-dualization functor. Then the following hold:
\begin{enumerate}
\item   $\BC\xra{[-,I]_r}\BC^\op$ is left adjoint to $\BC^\op\xra{[-,I]_l}\BC$; that is, the left and the right semi-dualization functors form a dual adjunction.
\item The  semi-dualization functors are normal monoidal functors $D_l,\ D_r\colon\BC^\op\rightarrow\BC^t$.
\end{enumerate}
\end{Proposition}
If $\BC$ is symmetric monoidal closed, one obviously has  $D:= [-,I]_l \simeq  [-,I]_r$ and, thus we only need to consider a single monoidal semi-dualization functor $D\colon\BC^\op\rightarrow\BC$.

Since monoidal functors send monoids to monoids, the identities   $\Comon\BC = \Mon(\BC^\op)^\op$ and $\Mon\BC = \Mon\BC^t$ imply
\begin{Corollary}\label{corr:}
Let $\BC$ be a  monoidal closed category. The semi-dualization functors induce functors 
 ${D_l},\ {D_r}\colon(\Comon \BC)^\op =\Mon\BC^\op\rightarrow\Mon\BC$, the \em left \em and the \em right dual monoid \em 
functor, such that the following diagram respectively commutes.
 \begin{equation}\label{diag}
\begin{minipage}{7cm}
\xymatrix@=2.5em{
\Mon\BC^\op\ar@/^0.5pc/@{->}[rr]^{{D_l}}\ar[d]_{|-|} \ar@/_0.5pc/@{->}[rr]_{D_r}
&&\Mon\BC\ar[d]^{|-|} \\ 
\BC^\op\ar@/^0.5pc/@{->}[rr]^{[-,I]_l }\ar@/_0.5pc/@{->}[rr]_{[-,I]_r }&&\BC 
}
\end{minipage}
\end{equation} 
\end{Corollary}

We call these functors \em dual monoid functors\em. 
For $\BC=\Mod_R$ for a commutative ring $R$ the functor  ${D}:= {D_l} = {D_r}$ is the so-called \em dual algebra	functor \em $$\Coalg_R^\op =\Comon(\Mod_R^\op)\rightarrow \Mon(\Mod_R)= \Alg_R,$$ while for
 $\BC={_R\Mod_R}$, the category of $R$-$R$-bimodules  over  a (not necessarily commutative) ring $R$,  the functors  ${D_l}$ and ${D_r}$ are the so-called \em dual ring	functors \em (see \cite{Swe}, \cite{LP_HP}) $$\Coring_R^\op =\Comon(_R\Mod_R^\op)\rightarrow \Mon(_R\Mod_R) = \ring_{R^\op}.$$
 
 \begin{Definition}\label{def:}\rm
An object $K$ in a symmetric monoidal closed category $\BC$ will be called {\em reflexive dualizable} provided that 
${D} K$ is a dual (in the monoidal sense) of $K$, by means of the evaluation $DK\otimes K\to I$, and that the adjunction counit {$K\xra{  \epsilon _K}DDK$} is an isomorphism. 
\end{Definition}
\begin{Remark}\label{rem:dualrd}\rm
The full subcategory $_{rd}\BC$ of $\BC$ spanned by all reflexive dualizable objects is a monoidal subcategory of $\BC$ and the monoidal semi-dualization functor $D$ induces, by restriction, a monoidal equivalence 
 $_{rd}\BC^\op\simeq  {_{rd}\BC}$.
The dual monoid functors lift this dual equivalence to an equivalence $\Mon{_{rd}\BC}^\op\simeq  \Mon{_{rd}\BC}$.
\end{Remark}

As is well known (see \cite{Sweedler} for example), for $\BC =\Vect$, the dual algebra functor $D\colon \Coalg_k^\op \rightarrow \Alg_k$ has a left adjoint  $(-)^\bullet\colon \Alg_k\rightarrow\Coalg_k^\op$, known as the  \em finite dual \em functor  or \em Sweedler dual \em functor.

Denoting, for an algebra $\ma$, the underlying vector space of the coalgebra $\ma^\bullet$ by $\ma^\circ$, the following hold:
\begin{enumerate}
\item $\ma^\circ$ is a subspace of $A^\ast$.
\item 
$\ma^\circ = A^\ast$ if $A$ is finite dimensional. 
In fact, more is true: The adjunction $(-)^\bullet\dashv D$ extends the dual equivalence ${_{fd}\Coalg_k}^\op\simeq {_{fd}\Alg_k}$ between the subcategories of finite dimensional algebras and coalgebras, respectively, lifted from the dual equivalence ${_{fd}\Vect}^\op\simeq {_{fd}\Vect}$.
\item If $\ma$ is the underlying algebra of a Hopf algebra then the coalgebra $\ma^\bullet$ carries a Hopf algebra structure.
\end{enumerate}

In fact, left adjoints of the dual monoid functor  exist more generally. The proof given for the case of $\BC =\Mod_R$, for $R$ a commutative ring (see \cite{HEP_AJSE},\cite{HEP_QMII}), can be generalized and so produces

\begin{Proposition}[\cite{HEP_QMI}]\label{prop:adjoints_l} %
The dual monoid functors of any monoidal closed  locally presentable category  $\BC$ have  left adjoints. 
\end{Proposition}

Since the categories $\Alg_R$  and  $\ring_R$ are locally finitely presentable (trivially), while the categories $\Coalg_R$ and 	$\Coring_R$ are locally  presentable (see \cite{HEP_QM},\cite{LP_HP}) (in fact locally countably  presentable by \cite{Ulmer}), we deduce: 

\begin{Corollary}\label{cor:lp}
\begin{enumerate}
\item For every commutative ring $R$, the  dual algebra  functor has a left adjoint $\mathsf{Sw} : \Alg_R \rightarrow \Coalg_R^\op$.
\item  For every ring $R$,
the left and right dual ring functors have left adjoints $\ring_{R^\op}\rightarrow \Coring_R^\op$.
\end{enumerate}
\end{Corollary}

We call these left adjoints \em generalized Sweedler dual functors\em.

\begin{Remark}\label{rem:}\rm
Note that, when $\BC$ is symmetric,  Proposition \ref{prop:adjoints_l} above is a special instance of the more general result, that in this case all so-called universal measuring comonoids exist and, thus, $\Mon\BC$ is enriched over $\Comon\BC$ (see \cite{Joy}, \cite{HV}).
\end{Remark}

\section{Left adjoints of induced functors}   
Let now  $G = (G,\Gamma,\gamma)\colon\BC\rightarrow\BD$    be a monoidal functor
 with  a left adjoint $L$ (which then is opmonoidal as  ${L} = (L,\Lambda,\lambda)$ by Remark \ref{rem:opmon})
 with unit $\sigma$  and counit $\rho$; we assume that 
its  induced functor ${G}\colon\Mon\BC\rightarrow\Mon\BD$ admits a left adjoint $F$ with unit $\eta$  and counit~$\epsilon$.

This situation is illustrated by the following diagram. 
 \begin{equation}\label{diag1}
\begin{minipage}{7cm}
\xymatrix@=2.5em{
\Mon\BC\ar[r]_{{G}}\ar[d]_{||-||}
&\Mon\BD\ar[d]^{|-|} \ar@/_0.8pc/@{->}[l]_{F} \\ 
\BC\ar[r]_{G}&\BD \ar@/_0.5pc/@{->}[l]_{L }
}
\end{minipage}
\end{equation}

We are, in this section,  concerned with the question, which of the properties of Sweedler's finite dual functor $F$ might share.

\subsection{The structure of a left adjoint}\label{ssec:struct}
Obviously, if such a left  adjoint exists, there is a 
 natural transformation $L|\sm|\xra{\kappa_\sm}||F\sm||$ characterized by commutativity  of the diagram
\begin{equation}\label{eqn:kappa}
\begin{minipage}{6cm}
\xymatrix@=2.5em{
{|\sm| }\ar[r]^{ \sigma_{|\sm|}} \ar[dr]_{|\eta_\sm|}& GL|\sm| \ar@{.>}[d]^{G\kappa_\sm }     \\
& { |{G}F\sm|}   = G||F\sm||  
}
\end{minipage}
\end{equation}

A straightforward calculation shows that, for $\BC = \Mod_R^\op$, $\BD = \Mod_R$ and $G$ the dualization functor (with $R=k$ a field), and $\ma = (A,m,e)$ any $k$-algebra, $\kappa_\ma$ is the embedding of the underlying vector space of the finite dual $\ma^\circ$ into the linear dual $A^\ast$.
Note however, that  there is no reason to believe that $\kappa_\ma$ is 
a monomorphism in $\Mod_R$ (that is, an epimorphism in $\BC$) in general.

In order to describe the algebraic properties 
 of $F\ma =(\ma^\circ, {m^\bullet}, e^\bullet)$, for some $\BD$-monoid $\ma = (A,m,e)$, we use the following concept.
 \begin{Definition}\label{def:indm}\rm
Let  $\ma = (A,m,e)$ be a $\BD$-monoid and $\bar{\ma}\ot \bar{\ma}\xra{ \bar{m}}\bar{\ma}$ and $I\xra{ \bar{e}} \bar{A}$
 $\BC$-morphisms such that the following diagram commutes for some $\psi$.
\begin{equation}\label{diag:SW4}
\begin{aligned}
\xymatrix@=2.5em{
   L(A\ot A)    \ar[rr]^{Lm}\ar[d]_{\Lambda_{A,A}}  &&{ LA} \ar[d]^{ \psi}  \\
LA\ot LA \ar[r]^{\psi\ot \psi}& 
\bar{\ma}\ot \bar{\ma}\ar[r]^{\bar{m}}  &\bar{\ma}
} 
\qquad 
\xymatrix@=2.5em{
LI  \ar[r]^{Le}\ar[d]_{\lambda}  &{ LA} \ar[d]^{ \psi}  \\
I \ar[r]^{\bar{e}}  &\bar{\ma}
} 
\end{aligned}
\end{equation}
Then the triple  $(\bar{\ma}, \bar{m}, \bar{e})$ is called  \em induced from $\ma$ \em by  $\psi$. 
If $(\bar{\ma}, \bar{m}, \bar{e})$ is  a monoid, then this monoid is called an {\em  induced quotient of  $\ma$ by $\psi$\em}\em\footnote{Note that $(\bar{\ma}, \bar{m}, \bar{e})$ is not necessarily a quotient of $\ma$ in the strict sense of the word. We will even use this term by abuse of language	if the morphism $LA\rightarrow A'$ is not epic.}.
\end{Definition}
Now the following holds.
\begin{Lemma}\label{lem:ind_quot}
If a left adjoint $F$ of the induced functor $G$ exists, then $F\ma$ is an induced quotient of $\ma$ by $\kappa_\ma$, for each $\BD$-monoid $\ma$.
\end{Lemma}

\begin{Proof}
The following diagram commutes where $\tilde{m}$ denotes the multiplication of $GF\ma$; the outer frame commutes since $\eta_\ma$ is a monoid morphism (use diagram \eqref{eqn:kappa}) and the small cells commute for obvious reasons.
\begin{equation}\label{diag:eta}
\begin{aligned}
\xymatrix@=3em{
   A  \ar[r]^{ \sigma_A }&{GLA } \ar[r]^{ G\kappa_\ma} & GF\ma \\
  A\ot A\ar[d]_{ \sigma_A\ot\sigma_A }   \ar[r]_{ \sigma_{A\ot A}} \ar[u]^{m}         &   GL(A\ot A)\ar[d]^{G\Lambda_{A,A}} \ar[u]^{GLm}& \\
  GLA\ot GLA\ar[r]^{\Gamma_{LA,LA}}\ar[rrd]_{G\kappa_\ma\ot G\kappa_\ma} & G(LA\ot LA)\ar[r]^{G(\kappa_\ma\ot\kappa_\ma)} & G(F\ma\ot F\ma)
\ar[uu]_{Gm^\bullet}\\
&&GF\ma\ot GF\ma\ar[u]^{\Gamma_{F\ma,F\ma}}\ar@/_3pc/@{->}[uuu]_{\tilde{m}}
}
\end{aligned}
\end{equation}
Thus, by the universal property of $\sigma_{A\ot A}$ one concludes that the left hand diagram in Definition \ref{def:indm} commutes, with {$m^\bullet$ replacing $\bar{m}$ and $\kappa_\ma$ replacing $\psi$.}

Concerning the right hand diagram observe first, that with $u$ the unit of $\si$ the outer frame as well as the left and top triangles of the following diagram commute, since $\eta_\si$ preserves units.
\begin{equation}\label{diag:}
\begin{aligned}
\xymatrix@=3em{
I \ar[r]^{ \gamma} \ar[dr]_{\sigma_I}\ar@/^1.5pc/@{->}[rr]^{\eta_\ma }  & GI\ar[r]^{Gu} & GF\si   \\
& GLI\ar[u]^{ G\lambda}   \ar[ur]_{G\kappa_\si}  &
}
\end{aligned}
\end{equation}
From this one similarly sees that the right hand diagram in Definition \ref{def:indm} commutes with the corresponding replacements.
\end{Proof}

\begin{Remark}\label{rem:indquot}\rm
\begin{enumerate}
\item Assume that $L$ is normal and $\psi$ is an epimorphism. {Suppose that the triple  $(\bar{\ma}, \bar{m}, \bar{e})$ is induced from $\ma$ by $\psi$ and that we know $\bar{m}$ is associative. Then $(\bar{\ma}, \bar{m}, \bar{e})$ is a monoid.  For, $\bar{e}$ is a left unit with respect to the multiplication $\bar{m}$ if the map $L(I\ot A)\xra{\Lambda_{I,A}}LI\ot LA\xra{\lambda\ot\psi}I\ot\bar{A}$ is an epimorphism; and similarly it is a right unit under the symmetric condition.}
\item   Every $\BC$-monoid $\mc$ is {an induced quotient of $G\mc$} by $\rho_C$, 
since  the following diagram 
 obviously commutes.
 \begin{equation}
 \begin{aligned}
\xymatrix@=2.5em{
  L(GC\ot GC) \ar[d]^{\Lambda_{GC,GC}}  \ar[r]^{L\Gamma_{C,C}  }&{ LG(C\ot C)}\ar[r]^{LGm}\ar[d]^{\rho_{C\ot C}} & LGC\ar[d]^{\rho_C }  \\
LGC\ot LGC  \ar[r]_{\rho_C\ot\rho_C}         &  C\ot C\ar[r]_m & C
}
\end{aligned}
\end{equation}

\item { Induced quotients of a monoid $\ma$ can be preordered as usual: $(LA \xra{ p} \mb) \leq (LA \xra{ q}\mc)$, if there is some $B \xra{ h}C$ with $h\circ p = q$. In particular one can talk about the \em smallest induced quotient of a monoid.\em}
\end{enumerate}
\end{Remark}

  We note for further use an equivalent description of induced quotients, available in module categories over a commutative ring (or, more generally, whenever the monoidal structure is given by universal bimorphisms), in the important case where $G$ is the semi-dualization functor.
 
For $R$-modules $M$ and $N$, the
unique map $ R^{M}\times R^{N}\xra{}R^{M\times N}$ making the following diagram commute, where $-\cdot -$ is the multiplication of $R$,

$$\begin{minipage}{6cm}
\xymatrix@=2em{M\times R^M\times N\times  R^N
     \ar@{.>}[r]^{ } \ar[d]_{ev_M\times ev_N  }&{ (M\times N)\times  R^{M\times N}} \ar[d]^{ ev_{M\times N} } \\
  R\times R  \ar[r]_{ -\cdot -}          &    R 
}
\end{minipage}$$
 is bilinear and so determines a linear map $\Pi_{M,N}\colon R^M\ot R^N\rightarrow R^{M\times N}$. 
 Analogously, one obtains a linear map $\Pi_{L,M,N}\colon R^L\ot R^M\ot R^N\rightarrow R^{L\times M\times N}$.
 With $t$ the universal bilinear map and $\theta$ the inclusion  $A^\ast\hookrightarrow R^A$ the following diagram (and a similar one for $\Pi_{A,A,A}$) commutes by definition of $\Lambda_{A,A}$ and $\Lambda_{A,A,A}$, respectively. 
  \begin{equation}\label{eqn:wis0}
 \begin{aligned}
\xymatrix@=2.5em{
A^\ast\ot A^\ast   \ar[rr]^{\theta_\ma\ot\theta_\ma}\ar[d]_{\Lambda_{A,A}} &&  R^A\ot R^A \ar[d]^{\Pi_{\ma,\ma}}\\
 (A\ot A)^\ast   \ar[r]^{\theta
 }
& R^{A\ot A}\ar[r]^{R^{t}}& R^{A\times A}
    }
       \end{aligned}
\end{equation}
Then, since $R^t\circ \theta$ is injective, in the diagram below the left hand upper cell commutes; that is, $m^\bullet$ is induced from $\ma$ by $\kappa_\ma$ if and only if the outer frame of the diagram commutes.
\begin{equation}\label{diag:wis1}
\begin{aligned}
\xymatrix@=3em{
  \ma^\circ   \ar[r]^{ m^\bullet }\ar[d]_{ \kappa_\ma }&{ \ma^\circ\ot\ma^\circ} \ar[r]^{ \kappa_\ma\ot  \kappa_\ma} & A^\ast\ot A^\ast\ar[rr]^{\theta_\ma\ot\theta_\ma}\ar@{.>}[d]^{\Lambda_{A,A}}&&R^A\ot R^A\ar[d]^{\Pi_{A,A}} \\
     A^\ast\ar[rr]^{m^\ast} \ar[d]^{\theta_\ma}& & (A\ot A)^\ast\ar[r]^{\theta_{A\ot A}}\ar[d]^{\theta_{A\ot A}} &  R^{A\ot A}\ar[r]^{R^t}   &     R^{A\times A} \\
     R^A \ar[rr]^{R^m} && R^{A\otimes A}\ar[rru]_{R^t} &&
}
\end{aligned}
\end{equation}

Motivated by the above we will use the following terminology {in which,} following Sweedler's notation,  we think of the functor $(-)^\circ$ as $\Mon\BD\xra{F}\Mon\BC\xra{||-||}\BC$.

\begin{Definition}\label{def:BS}\rm
We call the pair $((-)^\circ, \kappa)$ a \em basic situation \em where $(-)^\circ\colon\Mon\BD\rightarrow\BC$ is a functor and $\kappa\colon L\circ |-|\Rightarrow (-)^\circ$ is an epimorphic natural transformation\footnote{The requirement that $\kappa$ should be epimorphic is motivated by Example \ref{ex:triv} below.}.
\begin{equation*}
\xymatrix{
\Mon\BD \ar[rd]_{(-)^\circ}^(0.5){\phantom{a}}="1" \ar[rr]^{|-|}  && \BD \ar[ld]^{L}_(0.5){\phantom{a}}="2" \ar@{<=}"1";"2"^-{\kappa}
\\
& \BC 
}
\end{equation*}
\\
Any functor $\Mon\BD\xra{F}\Mon\BC$ with $||-||\circ F = (-)^\circ$ will be called a \em lift of $((-)^\circ, \kappa)$\em, provided that $F\ma$ is an induced quotient of $\ma$ by $\kappa_\ma$, for each $\BD$-monoid $\ma$.
\\
A basic situation admitting a lift will be called \em liftable (to $F$)\em.

When $L$ is the left adjoint of the semi-dualization functor $D$ of a symmetric monoidal closed category $\BD$, we call $(-)^\circ$ a \em Sweedler functor. \em 
\end{Definition}

Here is a simple example to show that a basic situation may fail to be liftable.
Let $G \colon\Mon\rightarrow\Set$ be the forgetful functor of the category of monoids. 
This is a normal symmetric  (cartesian) monoidal functor with the free monoid functor $F$ as left adjoint.
\begin{equation*}
\xymatrix{
& \Mon \ar[ld]^{F} \ar[rd]_{\Id}^(0.5){\phantom{a}}="1" \ar[rr]^{G}  && \BD \ar[ld]^{F}_(0.5){\phantom{a}}="2" \ar@{<=}"1";"2"^-{\epsilon}
\\
 \Mon_c \ar[rr]_{G} & & \Mon & 
}
\end{equation*}
With $(-)^\circ := \Id$ and $\kappa:=\epsilon$, the counit of the adjunction, we have a basic situation, where the other ${G}$ is the embedding of the category $\Mon_c$ of commutative monoids into $\Mon$. 
By the Eckmann-Hilton argument, $\sm^\circ$ can be supplied with a monoid structure only if $\sm$ is a commutative monoid.

The next lemma gives sufficient conditions for a lift. The notation
$\Lambda_{A,A,A}\colon L(A\ot A\ot A)\lra LA\ot LA\ot LA$ is used for the canonical map.
\begin{Lemma}\label{lem:Abullet}
Let $((-)^\circ,\kappa)$ be a basic situation. Assume  that, for each  $\BD$-monoid $\ma = (A,m,e)$, there exists a triple $\ma^\bullet:=(\ma^\circ, m^\bullet,e^\bullet)$ induced from $\ma$ by $\kappa_\ma$, and that the following conditions are satisfied.
\begin{enumerate}
\item The morphism $(\kappa_\ma\ot\kappa_\ma\ot \kappa_\ma)\circ\Lambda_{A,A,A}$ is an epimorphism.
\item The morphism 
$(\kappa_\ma\ot\kappa_\ma)\circ\Lambda_{A,A}$ is an epimorphism. 
\end{enumerate}
Then the following hold:
\begin{enumerate}
\renewcommand{\theenumi}{\alph{enumi}}
\item $\ma^\bullet$ is a monoid and {hence} an induced quotient of $\ma$ by $\kappa_\ma$, and this is the only induced quotient of $\ma$ by $\kappa_\ma$.
\item {The assignment $\ma\mapsto \ma^\bullet$ defines a lift of $((-)^\circ,\kappa)$.}
\end{enumerate}
\end{Lemma}
\begin{Proof}
Condition 2 implies by a diagram chase that, for any algebra {morphism} $\ma\xra{h}\mb$, the map $\ma^\circ\xra{h^\circ}\mb^\circ$ respects the comultiplication. Since $h^\circ$ preserves units trivially, this implies b.

Condition 1 implies coassociativity of $m^\bullet$ (the proof is the same {as in \cite[p. 113]{Swe}).}  That $e^\bullet$ is a counit, follows from item 1 of Remark \ref{rem:indquot}. 

The uniqueness statement is a consequence of the fact, that for any two lifts 
$(\ma^\circ, m^\bullet,e^\bullet)$ and $(\ma^\circ, m^\star,e^\star)$ induced from $\ma$ by $\kappa_\ma$ 
$\id_{\ma}^\circ$ would, by the above, {be a morphism} between them.
\end{Proof}

\subsection{The extension of dualities}\label{ssec:ext}

\subsubsection{The functor induced by a left adjoint}
It will prove useful to consider the following localized version of strong monoidality for a functor.
\begin{Definition}\label{def:strong}\rm
Let $G=(G,\Gamma,\gamma)\colon \BC\rightarrow \BD$ be a normal monoidal functor.
A $\BC$-object $C$ is called \em $G$-strong \em when $\Gamma_{C,C}$ and
$\Gamma_{C\otimes C,C}$ are invertible. 
The full subcategory of $\BC$ spanned by all $G$-strong objects will be denoted by $\BC_G$. A similar definition applies for a normal opmonoidal functor.
\end{Definition}

\begin{Proposition}\label{prop:normal}
Let $G = (G,\Gamma,\gamma)$ be a monoidal functor $\BC\rightarrow\BD$ admitting an adjunction  $F\dashv G$ with counit $\epsilon$.
Assume that the opmonoidal left adjoint  $F = (F,\Phi,\phi)$ of $G$ is normal.
Then the following hold for a $\BD$-monoid $\sfp = (P,p,u)$ whose underlying $\BD$-object $P$ is $F$-strong.
\begin{enumerate}
\item With $\bar{p} := FP\ot FP\xra{ \Phi_{P,P}^{-1}} F(P\ot P)\xra{ Fp} FP \text{\ \ and \ \ } \bar{u}:= I\xra{ \phi^{-1}}FI\xra{ Fu}FP$ the triple $\bar{\sfp} :=(FP, \bar{p},\bar{u})$ is a monoid in $\BC$. 

This defines a functor $\hat{F}\colon\Mon\BD_F\rightarrow\Mon\BC$\footnote{Concerning the notation $\Mon\BD_F$, we note that $\BD_F$ is not necessarily closed under the monoidal structure.}.
\item For every $\BC$-monoid $\sm = (M,m,e)$ the adjunction isomorphism $\BC(FP, M)\simeq \BD(P,GM)$ restricts to an isomorphism  $\Mon\BC(\hat{F}P,\sm)\simeq \Mon\BD(\sfp,{G}\sm).$
\end{enumerate}
\end{Proposition}

\begin{Proof}
The proof of the first statement is the same as for statement 2 in Remark \ref{rem:monfcts}.

To prove the second statement one needs to show that, for a monoid $\sfp$ in $\BD$ with $F$-strong $P$ and an arbitrary monoid $\sm$ in $\BC$, a map $f\colon P\rightarrow GM$ is a monoid morphism $\sfp\rightarrow {G}\sm$ if and only if the map $\bar{f}\colon FP\rightarrow M$ corresponding to $f$ under adjunction (that is, the map $FP\xra{ Ff}FGM\xra{ \epsilon_M}M$) is a monoid morphism $\bar{\sfp}\rightarrow \sm$.

In reference to the diagrams below, we thus see that commutativity of the left hand diagram is equivalent to that of the outer frame of the right hand diagram in both rows.

Since in the right hand diagram of (\ref{diag:hom1}) the upper left cell commutes by naturality of $\Phi$, the lower right cell by naturality of $\epsilon$, and the upper right cell by definition of $\Phi$ (see diagrams (\ref{eqn:mate2})), it is clear that the outer frame of this diagram commutes if the left hand diagram of (\ref{diag:hom1}) commutes.

Conversely, commutativity of the  right hand diagram implies the identity $$\epsilon_M\circ (F(Gm\circ \Gamma_{M,M}\circ (f\ot f))) = \epsilon\circ F(f\circ p).$$ Now commutativity of the left hand diagram follows by the universal	property of $\epsilon$.

\begin{equation}\label{diag:hom1}
\begin{aligned}
 \xymatrix{
P\ot P\ar[r]^{\!\!\! f\ot f}  \ar[dd]_{p  }&{ GM\ot GM} \ar[d]^{\Gamma_{M,M} } \\
             &   G(M\ot M)\ar[d]^{ Gm} \\ 
  P\ar[r]^f & GM
}
\qquad 
\xymatrix{
   FP\ot FP  \ar[r]^{\!\!\! Ff\ot Ff  }\ar[d]_{ \Phi_{P,P}^{-1} }&{FGM\ot FGM } \ar[r]^{\ \ \  \epsilon_M\ot \epsilon_M} & M\ot M\ar[ddd]^m \\
F(P\ot P)  \ar[r]_{F(f\ot f)\phantom{aa} }  \ar[dd]_{Fp}       &    F(GM\ot GM)\ar[d]_{F\Gamma_{M,M}}\ar[u]^{\Phi_{GM,GM}} &&\\
 & FG(M\ot M)\ar[d]^{FGm}\ar[uur]_{\epsilon_{M\ot M}} &\\
FP\ar[r]^{Ff} & FGM\ar[r]^{\epsilon_M} & M
}
\end{aligned}
\end{equation}

\begin{equation}\label{diag:hom2}
\begin{aligned}
\xymatrix@=3em{
I\ar[r]^u\ar[d]_\gamma & P\ar[d]^f\\
GI\ar[r]^{Ge} & GM
}\qquad 
\xymatrix@=3em{
I\ar[dd]_{id}\ar[r]^{\phi^{-1}} & FI\ar[d]^{F\gamma}\ar[r]^{Fu} & P\ar[d]^{Ff}\\
& FGI\ar[r]^{FGe}\ar[dl]^{\epsilon_I} & FGM\ar[d]^{\epsilon_M}\\
I\ar[rr]^e && M
}
\end{aligned}
\end{equation}
The proof for unit preservation is similar.
\end{Proof}

\begin{Remark}\label{rem:crux}\rm
\begin{enumerate}
\item Let $\sfp$ be a $\BD$-monoid with $F$-strong $P$. 
Since the unit $\eta_P$ corresponds under adjunction to $id_{FP}$, 
it follows from Proposition \ref{prop:normal} that $\eta_P$ is a monoid morphism 
$\sfp\rightarrow {G}\bar{\sfp}$ and, in fact, a ${G}$-universal morphism for $\sfp$.
\item If the opmonoidal left adjoint $F$ of $G$ is strong and therefore can be seen as a monoidal functor as well, statement 2 above gives the familiar		result that,
given a monoidal adjunction $F \dashv G$ with unit $\eta\colon\Id\rightarrow GF$ and counit $\epsilon\colon FG\rightarrow\Id$,  the induced functors ${F}$ and ${G}$ form an adjunction with unit $\eta'$ and counit $\epsilon'$, such that   $|\eta' | = \eta $ and $|\epsilon'| = \epsilon$. 
\end{enumerate}
\end{Remark}

\subsubsection{Induced dualities}\label{dual_ind}

Every  adjunction  $\BA\xrightarrow{G}\BB\xrightarrow{F}\BA$  with unit $\eta:\Id_\BB\Rightarrow GF$ and counit $\epsilon\colon FG\Rightarrow\Id_\BA$  induces, by restriction and corestriction,  an equivalence $Fix\eta\simeq Fix\epsilon$ between the full and replete subcategories of $\BB$ and $\BA$, respectively, spanned by those objects whose components of $\eta$ and $\epsilon$ respectively are isomorphisms. 

When dealing with a monoidal functor ${\BA}\xra{ G}{\BB}$ 
it seems natural to relate this equivalence to the equivalence induced by the adjunction $F\dashv {G}$, provided that the latter exists.

As a paradigmatic example consider the adjunction for the (monoidal) dualization functor $(-)^\ast\colon\Mod_R^\op\rightarrow\Mod_R$, which restricts to a (dual) equivalence ${_{f\!gp}\Mod_R}^\op\simeq {_{f\!gp}\Mod_R}$ for 
 the category   ${_{f\!gp}\Mod_R}$ of finitely generated projective $R$-modules. This equivalence can be lifted to the monoid level and yields the duality ${_{f\!gp}\Coalg_R}^\op\simeq {_{f\!gp}\Alg_R}$ for the categories  of finite dimensional $R$-algebras and  coalgebras respectively. 
 Since the functor ${\Coalg_R}^\op\rightarrow {\Alg_R}$ induced by $(-)^\ast$ has a left adjoint (see Corollary \ref{cor:lp}), one may ask whether the (dual) equivalence ${_{f\!gp}\Coalg_R}^\op\simeq {_{f\!gp}\Alg_R}$ is induced by this adjunction. 
 If  $R=k$ is a field this is known to be true (see  above). 
We now will show that this question can be answered in the affirmative for every commutative ring $R$.

By Remark \ref{rem:monfcts}, one first obtains:

\begin{Lemma}\label{lem:1}
For $A\in Fix\epsilon\cap\BA_G$ one has  {$GA\in \BB_F$. Consequently, 
$GA\in Fix\eta\cap\BB_F$ if and only if $A\ot A\in Fix\epsilon$.}
\end{Lemma}

As a consequence, using notation as in Remarks \ref{rem:monfcts} and Proposition \ref{prop:normal}, we obtain:
\begin{Proposition}\label{prop:restr_fct}
Let $\mathbf{A}\xra{ G}\mathbf{B}$ be a monoidal functor with normal opmonoidal left adjoint $L$.
Assume that the induced functor $G$ has a left adjoint $F$.

Consider any full and replete subcategory $\BA'$ of $Fix\epsilon\cap\BA_G$, closed under tensor squaring, and denote by $\BB'$ the full and replete subcategory of $\BB$ spanned by the objects $GA$ for $A\in \BA'$.
Then
\begin{enumerate}
\item  the functors $\BA \xra{ G} \BB$ and $\BB \xra{ L}\BA$   
can be restricted to functors $\BA'\xra{ G'}\BB'\xra{ L'}\BA'$, and these provide an equivalence $\BA'\simeq\BB'$.
\item  the functors $\Mon\BA \xra{ G} \Mon\BB $ and $\Mon\BB_L\xra{ \widehat{L}}\Mon\BA$ can be restricted to functors 
$\Mon\BA'\xra{ {G}'}\Mon\BB'\xra{\widehat{L}'}\Mon\BA'$, and these  
provide an equivalence $\Mon\BA'\simeq\Mon\BB'$. 
\item $ \widehat{L}'$ is the restriction of $F\colon\Mon\BB\lra\Mon\BA $ to the respective subcategories.
\end{enumerate}
\end{Proposition}

\begin{Proof}
By the lemma above $G$ can be restricted to a functor $\BA'\rightarrow Fix\eta\cap\BB_F$; thus, $\BB'$ is a subcategory of $Fix\eta\cap\BB_F$. If $B$ belongs to $\BB'$, one has $B\simeq GA$ for some $A$ in $\BA'$ and thus $B\ot B\simeq GA\ot GA\simeq G(A\ot A)$ (because $A$ is $G$-strong). 
Hence one can apply the lemma to $L$ as well and item 1 follows.

$\widehat{L}$ as defined in Proposition \ref{prop:normal} exists such that  
  $\widehat{L}'$ is well defined as the restriction and corestriction. 
  It follows from Remark \ref{rem:crux}.1 that, for each $B\in\BB'$, the unit $\eta_B$ lifts to an  isomorphism $B\rightarrow {G}'\widehat{L}'B={G}\widehat{L}B$ which, moreover, is ${G}$-universal. Since  ${G}'$ can alternatively be seen as obtained by restriction and corestriction of the functor $\widehat{G}$ in the sense of Proposition \ref{prop:normal}, it follows the same way that the counit $\epsilon$ lifts to a natural isomorphism $ \widehat{L}'{G}'\rightarrow id$. This proves item 2.  
  
   Finally, item 3 is a consequence of Proposition \ref{prop:normal} item 2.
\end{Proof}

Applying this to the left semi-dualization functor $G = [-,I]_l$ in interesting cases we obtain: 
\begin{Examples}\label{fact:exs}\rm
\begin{enumerate}
\item Let $\BC$ be a symmetric monoidal closed subcategory. By Remark \ref{rem:dualrd} 
the  subcategory $_{rd}\BC$ and its dual can be chosen as $\BB'$ and $\BA'$ respectively in the proposition above.

\item Let $R$ be a commutative ring. Then $G$ is the dualization functor $(-)^\ast\colon\Mod_R^\op\rightarrow\Mod_R$. 
Again by  Remark \ref{rem:dualrd} we have 
 ${_{f\!gp}\Mod_R}= {_{rd}\Mod_R}$ and, since  the dual of a finitely generated projective module is finitely generated projective again, ${_{f\!gp}\Mod_R}$ and its dual can be chosen as $\BB'$ and $\BA'$, respectively, in the proposition above.

Consequently, Proposition \ref{prop:restr_fct} yields the familiar duality ${_{f\!gp}\Coalg_R}^\op\simeq {_{f\!gp}\Alg_R}$ mentioned above.
The  duality ${_{f\!d}\Coalg_k}^\op\simeq {_{f\!d}\Alg_k}$ for the categories of finite dimensional $k$-algebras and coalgebras respectively over a field $k$ is a special instance. 

\item   Let $R$ be a not necessarily commutative ring. Then $G$ is the left	semi-dualization functor
for the category $_R\Mod_R$ of $R$-$R$-bimodules  (see  \cite{LP_HP}), which has the right semi-dualization functor as a left adjoint (see  \cite{LP_HP}). Similarly as in the commutative case  the following holds for 
 ${_R\mathbf{FGP}}$ and $\mathbf{FGP}_R$, the categories  of bimodules which are finitely generated projective if considered as left and right $R$-modules, respectively: 
$G$ maps ${_R\mathbf{FGP}}^\op$ into $\mathbf{FGP}_R$, and
${_R\mathbf{FGP}}^\op \subset Fix\epsilon\cap {_R\Mod_R}_G$ and  $\mathbf{FGP}_R \subset Fix\eta\cap {_R\Mod_R}_F$.
 
 The duality ${_R\mathbf{FGP}}^\op\simeq \mathbf{FGP}_R$ resulting by  Proposition \ref{prop:restr_fct} is the duality 
 $_{f\!gp}\Coring_R^\op$ $\simeq {_{f\!gp}\ring_{R^\op}}$
 for the categories  of $R$-corings and $R^\op$-rings,  which are finitely generated projective as left $R$-modules and  as right $R$-modules, respectively. This duality coincides with Sweedler's duality (\cite{Swe}) given by the dual $R$-ring functor and the dual $R$-coring construction for finitely generated projective $R$-rings (see \cite{LP_HP}).
\end{enumerate}
\end{Examples}

\begin{Remark}\label{rem:inddual}\rm
These examples show in particular, that the question asked at the beginning of this section, whether for each commutative ring $R$ the familiar duality ${_{f\!gp}\Coalg_R}^\op\simeq {_{f\!gp}\Alg_R}$ is, as in the field case, a restriction	of the adjunction $\mathsf{Sw} \dashv{G}$ can be answered in the affirmative. And similarly for the non-commutative case and the dual coring adjunction.

 If one considers this a property	 to be shared by any reasonable generalization of the Sweedler dual functor,  
 it seems justified to call the left adjoints of the dualization functors, 
 which exist by Proposition~\ref{prop:adjoints_l},  {\em generalized Sweedler duals}, as we did above. 
\end{Remark}

\subsection{The lift to bimonoids}\label{subsec:com}

\subsubsection{Symmetric monoidal functors}\label{ssec:sym}
Recall that 
a monoidal functor $F$ between symmetric monoidal categories is called \em symmetric \em if the following diagram commutes.
$$\begin{minipage}{6cm}
\xymatrix@=2em{
    FC\ot FD \ar[r]^{\Phi_{C,D}  }\ar[d]_{ \sigma_{} }&{ F(C\ot D)} \ar[d]^{F\sigma }  \\
  FD\ot FC  \ar[r]_{ \Phi_{D,C}}          &    F(D\ot C)
}
\end{minipage}$$
\em Symmetric opmonoidal \em is defined dually.

If $(G,\Gamma,\gamma)$ is symmetric monoidal then its opmonoidal left adjoint $(F,\Phi,\phi)$, if it exists, is symmetric.

For every commutative ring $R$ the dualization functor $(-)^\ast$ is a symmetric monoidal functor as a simple calculation shows. 
This can be generalized as we will now show. 

\begin{Lemma}\label{lem:swing}
For  $f\colon X\ot A\rightarrow I$ and $B,Y$ in any symmetric monoidal category, the following diagram commutes.
$$\begin{minipage}{6cm}
\xymatrix@=2.5em{
 X\ot Y\ot A\ot B    \ar[rr]^{ id\ot \sigma^{-1}_{A,Y}\ot id }\ar[d]_{\sigma_{X,Y}\ot id\ot id} & &{ X\ot A\ot Y\ot B} \ar[d]^{f\ot id\ot id}  \\
  Y\ot X\ot A\ot B     \ar[rr]_{ id \ot f\ot id}    &      &    Y \ot B
}
\end{minipage}$$
\end{Lemma}

\begin{Proof}
Interpreting the square as a string diagram: 
 \begin{center}
\psscalebox{0.7 0.7} 
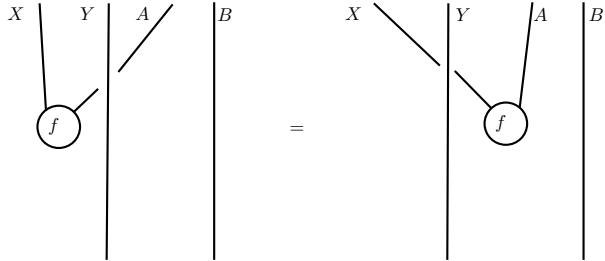
{
\begin{pspicture}(0,-2.4512093)(11.56,2.4512093)
\pscircle[linecolor=black, linewidth=0.04, dimen=outer](0.96,0.08882196){0.42}
\rput[bl](0.76,-0.051178034){$\Huge{f}$}
\pscircle[linecolor=black, linewidth=0.04, dimen=outer](9.36,0.14882196){0.42}
\rput[bl](9.16,0.008821964){$\Huge{f}$}
\rput[bl](5.3,-0.031178035){$\Huge{=}$}
\psline[linecolor=black, linewidth=0.04](0.6,2.428822)(0.72,0.38882196)
\psline[linecolor=black, linewidth=0.04](1.9,2.428822)(1.86,-2.3911781)(1.86,-2.431178)
\psline[linecolor=black, linewidth=0.04](3.1,2.408822)(2.08,1.128822)
\psline[linecolor=black, linewidth=0.04](1.24,0.36882198)(1.7,0.808822)
\psline[linecolor=black, linewidth=0.04](3.88,2.448822)(3.88,-2.431178)
\psline[linecolor=black, linewidth=0.04](8.28,2.428822)(8.26,-2.431178)(8.26,-2.431178)
\psline[linecolor=black, linewidth=0.04](9.62,0.44882196)(9.86,2.448822)
\psline[linecolor=black, linewidth=0.04](6.88,2.428822)(8.12,1.248822)
\psline[linecolor=black, linewidth=0.04](8.4,1.148822)(9.08,0.44882196)
\psline[linecolor=black, linewidth=0.04](10.82,2.448822)(10.82,-2.431178)
\rput[bl](6.34,2.1088219){$X$}
\rput[bl](8.42,2.088822){$Y$}
\rput[bl](9.88,2.088822){$A$}
\rput[bl](10.92,2.088822){$B$}
\rput[bl](0.0,2.1088219){$X$}
\rput[bl](1.36,2.1088219){$Y$}
\rput[bl](2.4,2.1088219){$A$}
\rput[bl](3.94,2.088822){$B$}
\end{pspicture}
}
\end{center}
provides immediate proof.
\end{Proof}

The following now is an immediate consequence of Proposition \ref{prop:dual-mon} item 2.

\begin{Proposition}\label{prop:dualsym}
The semi-dualization functor of a symmetric monoidal closed category is a normal symmetric monoidal functor.
\end{Proposition}

Let $(G,\Gamma,\gamma)$ be a symmetric monoidal functor.
By definition $\Gamma$ is a natural transformation $G_1\Rightarrow G_2$ between	the functors 
$G_1=\BC\times\BC\xra{G\times G}\BD\times\BD\xra{-\otimes -}\BD \text{\ and \  }  G_2=\BC\times\BC\xra{-\otimes -}\BC\xra{G}\BD $
 which are symmetric monoidal, since symmetric monoidal functors compose.

For every commutative ring $R$ the symmetric {monoidal semi-dualization functor} $((-)^\ast, \Lambda,\lambda)$ has the property that $\Lambda\colon (-)^\ast_1\Rightarrow (-)^\ast_2$ is a monoidal transformation, 	
as a simple calculation shows. This can be generalized as follows:
\begin{Lemma}\label{lem:sym_mon_trans}
If  $G=(G,\Gamma,\gamma)\colon\BC\rightarrow\BD$ is a symmetric monoidal functor then $\Gamma\colon G_1\Rightarrow G_2$ is a monoidal transformation.
\end{Lemma} 

\begin{Proof}
We may assume (see \cite{Day}) that $\BD$ is monoidally cocomplete. Thus,   monoidal functors  correspond to monoids and monoidal transformations to monoid morphisms (see Remark~\ref{rem:day}). By this correspondence  symmetric monoidal functors correspond to commutative monoids.
Since the multiplication of a commutative monoid is a monoid morphism, the claim follows.
\end{Proof}

Recall that, for $\BC$ symmetric monoidal, the tensor product on $\Mon\BC$ is the functor induced by the monoidal functor  $\BC\times\BC\xra{-\otimes -}\BC$. One then has, for a symmetric monoidal functor $(G,\Gamma,\gamma)\colon \BC\rightarrow\BD$,  induced functors as follows.
\begin{eqnarray*}
{G}_1 &=& \Mon\BC\times\Mon\BC\xra{{G}\times {G}}\Mon\BD\times\Mon\BD\xra{-\otimes -}\Mon\BD\\
{G}_2 &=& \Mon\BC\times\Mon\BC\xra{-\otimes -}\Mon\BC\xra{{G}}\Mon\BD
\end{eqnarray*}

By Remark \ref{rem:monfcts} item 4, we thus obtain the following as a corollary (also see \cite[Lemma 1.4.2]{Schau} where, however, no proof is given).

\begin{Theorem}\label{prop:ind_fct_mon}
Let  $G = (G,\Gamma,\gamma)\colon\BC\rightarrow\BD $  be a symmetric monoidal functor. Then,  for each pair of monoids
 $\sm$ and $\sn$  in $\BC$, the map $\Gamma_{M,N}$ is a monoid morphism ${G}\sm\otimes   {G}\sn \rightarrow {G}(\sm   \otimes\sn)$. Moreover, $\gamma$ is a monoid morphism $\si\rightarrow{G}\si$.
This makes the induced functor ${G}\colon\Mon\BC\rightarrow\Mon\BD$ a symmetric monoidal functor
with multiplications	 $\Gamma_{\sm,\sn} = \Gamma_{M,N}$ and unit $\gamma$. 
\end{Theorem}

As a corollary we obtain:
\begin{Proposition}\label{con:con1}
The  dual monoid functor  
 ${D}\colon(\Comon\BC)^\op\rightarrow\Mon\BC$ of a symmetric monoi\-dal category $\BC$ is a symmetric monoidal functor.
\end{Proposition}

\begin{Remarks}\label{rem:Xi}\rm
\begin{enumerate}
\item
If the functor ${G}\colon\Mon\BC\rightarrow\Mon\BD$ has a left adjoint $F$, this  is an opmonoidal functor with a natural transformation $\Phi_{\ma,\mb}\colon F(\ma\ot \mb)\rightarrow F\ma\ot F\mb$ making the following diagram commute. (The left upper cell commutes by definition of $\Lambda$ and the bottom cell commutes by naturality of $\Gamma$. The outer frame commutes by the definition of mates.)

$$\begin{minipage}{6cm}
\xymatrix@=2.5em{
  A\ot B   \ar[r]^{ \lambda_{A\ot B} }\ar[d]_{\lambda_A\ot\lambda_B  }&{ GL(A\ot B)} \ar[d]^{ G\Lambda_{A,B}} \ar[r]^{G\kappa_{\ma\ot \mb}} & GF(\ma\ot \mb)\ar[dd]^{G\Phi_{\ma,\mb}}\\
GLA\ot GL B   \ar[r]_{ \Gamma_{LA,LB}} \ar[d]_{G\kappa_\ma\ot G\kappa_\mb}         &    G(LA\ot LB) \ar[dr]^{G(\kappa_{\ma}\ot\kappa_\mb)}\\
GF\ma\ot GF\mb\ar[rr]^{\Gamma_{F\ma , F\mb}} && G(F\ma\ot F\mb)
}
\end{minipage}$$
The universal property of	$\lambda_{A\ot B}$ now implies commutativity of the diagram
\begin{equation*}
\begin{minipage}{6cm}
\xymatrix@=2.5em{
   LA\ot LB  \ar[d]_{\kappa_{\ma}\ot\kappa_{\mb}  }&  L(A\ot B) \ar[l]_{\ \Lambda_{A,B} } \ar[d]^{ \kappa_{\ma\ot\mb}}  \\
 {F\ma\ot F\mb }          &   F(\ma\ot\mb)  \ar[l]^{\ \ \Phi_{\ma,\mb}}  
}
\end{minipage}
\end{equation*}
and $\Phi_{\ma,\mb}$ is (under the given hypotheses) the only natural transformation	making this diagram commute. 
\item If $G$ is the dual monoid functor on $\BD$ and the generalized Sweedler dual $\mathsf{Sw}$ exists
then the above is a weakening in the sense of (b) of \cite[Lemma 6.0.1]{Sweedler} 
in that we only get a morphism	$\mathsf{Sw} A\ot \mathsf{Sw} B\lra  \mathsf{Sw}(A\ot B)$ in $\BD$, but not necessarily an isomorphism. 
 As shown above, this requirement is not necessary for obtaining the desired left adjoint. 
 Item (b) of \cite[Lemma 6.0.1]{Sweedler} holds if and only if $\mathsf{Sw} $, considered  as the opmonoidal left adjoint of ${G}$, is strong opmonoidal.
\end{enumerate}
\end{Remarks}

\subsubsection{The generalized Sweedler dual of a bimonoid}\label{ssec:genfindualbim}

Recall that, for  a symmetric monoidal category $\BC$, one has  the category  of bimonoids in $\BC$,  defined as $\Bimon\BC:= \Mon(\Comon\BC)  = \Comon(\Mon\BC)$. 
For bimonoids we use the notation $\mathsf{B} = (\mathsf{B}^a,\mathsf{B}^c)$, 
where the components denote the underlying monoid and comonoid, respectively.

Since opmonoidal functors send comonoids to comonoids one obtains
\begin{Proposition}\label{prop:bimon}
In the situtation depicted in diagram (\ref{diag1}), the opmonoidal left adjoint $F$ of the induced functor ${G}$ maps bimonoids to bimonoids; that is, $F$ induces a functor
${F}\colon\Bimon\BD\rightarrow\Bimon\BC$ such that the following diagram commutes.
$$\begin{minipage}{6cm}
\xymatrix@=2.5em{
 \Bimon\BD    \ar[r]^{ {F} }\ar[d]_{|-|  }&{\Bimon\BC } \ar[d]^{|-| }  \\
  \Mon\BD   \ar[r]_{ F}          &   \Mon\BC
}
\end{minipage}$$
 \end{Proposition}
  
 Specializing this to the dualization functor $\Mod_R^\op\rightarrow\Mod_R$, in view of Proposition~\ref{con:con1} and $\Bialg_R := \Bimon(\Mod_R)$, we deduce:

\begin{Proposition}\label{corr:bialg}
Let $R$ be a commutative ring. Then the generalized Sweedler dual functor $\mathsf{Sw} \colon\Alg_R\rightarrow\Coalg_R^\op$ maps bialgebras to bialgebras.
\end{Proposition}

As mentioned before, Sweedler's  finite dual functor sends Hopf algebras to Hopf algebras. The following is a generalization.

\begin{Proposition}\label{prop:hopf}
If, in the situation of Proposition \ref{prop:bimon}, a bimonoid $\mh =(\mh^a,\mh^c)$ is a Hopf monoid with antipode $S\colon\mh^a\rightarrow(\mh^a)^\op$, then ${F}\mh$ is a Hopf monoid with antipode $FS$, provided that $\kappa_{\mh^a}$ is an epimorphism in $\BC$.
\end{Proposition}

\begin{Proof}
Let $\mh = (\mh^a,\mh^c,S)$ be a Hopf monoid over $\BD$, where $\mh^a = (H,m,e)$ is the underlying monoid of $\mh$ and $\mh^c = (H,{\mu},\epsilon)$ is its underlying comonoid. 
\\
Put
${F}\mh = (F\mh^a,\tilde{\mu},\tilde{\epsilon})$ with
\begin{enumerate}
\item $F\mh^a =: (C,c,n)$
\item $\tilde{\mu} = C\xra{ F\mu}F(\mh^a\otimes\mh^a)\xra{ \Phi_{\mh^a,\mh^a}}C\otimes C$
\item $\tilde{\epsilon} = C\xra{Fe}FI\xra{ \phi} I$
\end{enumerate}
Then ${G}(F\mh^a) = (GC, \tilde{m},\tilde{n})$ with $\tilde{m} = GC\otimes GC\xra{ \Gamma_{C,C}} G(C\otimes C)\xra{ Gc} GC$ and \\Ê $\tilde{n} = I\xra{\gamma }GI\xra{ Gn}GC$.

Since the antipode $S$ of $\mh$ is a monoid morphism  $\mh^a\rightarrow  (\mh^a)^\op$ in $\BD$, there is the monoid morphism $FS\colon F\mh\rightarrow F(\mh^a)^\op$.  To show this map is an antipode for ${F}\mh$, we must check the identity:
\begin{eqnarray*}
C\xra{ F\mu}F(\mh^a\otimes\mh^a)\xra{ \Phi_{\mh^a,\mh^a}}C\otimes C\xra{{FS}\ot id}C\ot C\xra{c} C &=& C\xra{Fe}FI\xra{\phi}I\xra{n }C\\
\end{eqnarray*}
Using the definitions above and that of an opmonoidal left adjoint, and also the facts that $\eta_\mh$ is a monoid morphism and $S$ satisfies the antipode equations for $H$, one sees that the following diagram commutes.

$$\begin{minipage}{6cm}\small
\xymatrix@=3em{
  GC   \ar[r]^{\!\!\!\! \!\!\!\! \!\!\!\! \!\!\!\! GF\mu }&{ GF(\mh^a\otimes\mh^a)} \ar[r]^{ G\Phi_{\mh^a,\mh^a} \!\!\!\! \!\!\!\!} & G(C\otimes C)\ar[r]^{G(\tilde{S}\ot id)} & G(C\otimes C)\ar[r]^{Gc} &GC \\
  & & GC\ot GC \ar[r]_{GFS\ot id}\ar[u]_{\Gamma_{C,C}} &GC\ot GC\ar[u]_{\Gamma_{C,C}} \ar[ur]_{\tilde{m}} &\\
   &    H\ot H \ar[r]_{S\ot id} \ar[ur]_{\eta_{\mh^a}\ot \eta_{\mh^a}}\ar[uu]^{\eta_{\mh^a\ot\mh^a}}  &   H\ot H\ar[ur]_{\eta_{\mh^a}\ot\eta_{\mh^a}} \ar[drr]^m \\
 H\ar[uuu]^{\eta_{\mh^a}}  \ar[d]_{\eta_{\mh^a}}   \ar[rr]_{ \epsilon} \ar[ur]^\mu        & &    I\ar[rr]_e\ar[d]^{\eta_{\mathsf{I}}}\ar[dr]^\gamma & & H \ar[d]^{\eta_{\mh^a}} \ar[uuu]_{\eta_{\mh^a}} \\
 GC    \ar[rr]^{ GFe } &&{GFI } \ar[r]^{G\phi } &GI\ar[r]^{Gn} & GC
}
\end{minipage}$$
To prove the equality of the two maps $FH=C\to C=FH$, since 
$\kappa_{\mh^a} : LH\to FH$ (see Diagram (5))
is an epimorphism, it the same as proving two maps $LH \to FH$
are equal. By adjointness, this is the same as proving two maps
$H \to GFH$ are equal. That is what the above diagram does.
%
\end{Proof}

As applications we mention the following examples.
\begin{Examples}\label{ex:triv}\rm
\begin{enumerate}
\item For every symmetric monoidal functor $G$  with a left adjoint, the opmonoidal left adjoint of its induced functor ${G}$, if it exists, maps those Hopf algebras $(\mh,S)$, for which the canonical map $\kappa_{\mh^a}$ is injective, to Hopf algebras. 

In particular, for every commutative ring $R$, the generalized finite dual functor $\mathsf{Sw}$ maps those Hopf algebras $(\mh,S)$, for which the canonical map $\kappa_{\mh^a}$ is injective, to Hopf algebras. This implies that the original finite dual functor (over a field) maps Hopf algebras to Hopf algebras.

\item If $L\dashv G$ is  a monoidal equivalence then the opmonoidal left adjoint of ${G}$   maps Hopf algebras $(\mh,S)$ to Hopf algebras. 
In particular, we obtain the familiar fact that, by dualization, one can assign a Hopf algebra to any Hopf algebra whose underlying module is finitely generated projective.
\end{enumerate}
\end{Examples}

 We end this section with the observation that in quite a number of instances the dual bimonoid functor  of Proposition \ref{prop:bimon} 
is again part of a dual adjunction. In fact, since the forgetful functors create colimits (for example, see \cite{HEP_QM}) and $F\colon\Mon\BD\rightarrow\Mon\BC$, being a left adjoint, preserves these, the functor $F\colon\Bimon\BD\rightarrow \Bimon\BC$ also preserves colimits.
If now the category $\BD$ is symmetric monoidal closed and locally presentable then $\Bimon\BD$ is locally presentable again (see again\cite{HEP_QM}) and, thus, $F$ has a right adjoint by the Special Adjoint Functor Theorem. 
If $F$ moreover maps Hopf monoids to Hopf monoids, the same argument applies to 
the functor $F\colon\Hopf\BC\rightarrow \Hopf\BD$ between the categories of Hopf 
monoids in $\BC$ and $\BD$ respectively, provided that, in addition, tensor squaring in these categories preserves extremal epimorphisms, since then these categories are coreflective in their categories of bimonoids (see \cite[Thm. 54]{HEP_QMI}). 
We thus have obtained:
 \begin{Proposition}\label{prop:dualadj2}
Let $G\colon\BC\rightarrow\BD$ be a symmetric monoidal functor, where $\BC$ and $\BD$ are symmetric monoidal closed and $\BD$ is locally presentable. Let $F\colon\Mon\BD\rightarrow\Mon\BC$ be the opmonoidal left adjoint of the functor induced by $G$. 
Then:
\begin{enumerate}
\item the functor $F$, considered as a functor $\Bimon\BD\rightarrow\Bimon\BC$, has a right adjoint;
\item if the natural transformation $\kappa$ is epimorphic and tensor squaring in $\BC$ 
and $\BD$ preserves extremal epimorphisms then the functor $F$, 
considered as a functor $\Hopf\BD\rightarrow\Hopf\BC$, has a right adjoint.\end{enumerate}
\end{Proposition}

\section{Constructing a left adjoint}\label{sec:con}

\subsection{A criterion for adjointness}\label{ssec:adj}
 
\begin{Theorem}\label{prop:leftadjointgeneral}
Given a basic situation $((-)^\circ,\kappa)$ (see Definition~\ref{def:BS}), the functor ${G}\colon \Mon\BC\lra\Mon\BD$ has a left adjoint $F $, 
provided that  
 
\begin{enumerate}
\item $((-)^\circ,\kappa)$ is liftable (to a functor $(-)^\bullet$).
\item $\ma^\bullet$ is the smallest induced quotient of $\ma$, for each  $\BD$-monoid $\ma$,
\item the morphisms $\kappa_{\ma}\ot\kappa_{ {\ma}}$ are epimorphisms \footnote{Note that this condition is satisfied if $\BC$ is monoidal closed.}.
\end{enumerate}
Then $F\ma = \ma^\bullet$ for each $\BD$-monoid $\ma$.
\end{Theorem}

\begin{Proof}
By item 1 there exists a functor $(-)^\bullet\colon \Mon\BD \rightarrow\Mon\BC$.
That this functor is left adjoint to ${G}$ is a consequence of the following lemmas.
\end{Proof}

\begin{Lemma}\label{lem:2}
For every $\BC$-monoid $\mc = (C,m,e)$, the counit $\rho_C$ factors 
as $LGC\xra{\kappa_{{G}\mc}}{({G}\mc)^\circ}\xra{ {\epsilon}_\mc}~\!\!C$ with a monoid morphism  ${\epsilon}_\mc\colon ({G}\mc)^\bullet\rightarrow \mc$.
\end{Lemma}
\begin{Proof}
With $\ma = {G}\mc$ diagram (\ref{diag:SW4}) becomes diagram \eqref{diag:LGC}. 
\begin{equation}\label{diag:LGC}
\begin{aligned}
\xymatrix@=2.5em{
  L(GC\ot GC) \ar[d]^{\Lambda_{GC,GC}}  \ar[r]^{L\Gamma_{C,C}  }&{ LG(C\ot C)}\ar[r]^{LGm} & LGC\ar[d]^{\kappa_{{G}\mc} }  \\
LGC\ot LGC  \ar[r]_{\kappa_{{G}\mc} \ot\kappa_{{G}{\mc}}}         &  ({G}\mc)^\circ\ot ({G}\mc)^\circ \ar[r]^{\ \ \ m_{{G}\mc}} & ({G}\mc)^\circ 
}
\end{aligned}
\end{equation}
Since $\mc$ is induced by $\rho_C$ (see item 1 of Remark \ref{rem:indquot}), 
we conclude from the {minimality condition 1(b) on $(G\mc)^\bullet$} 
that there exists a unique morphism $\epsilon_\mc\colon  ({G}\mc)^\bullet \rightarrow C$ with
$\epsilon_\mc \circ \kappa_{{G}\mc} = \rho_C \text{\ \ and\ \ } \epsilon_\mc \circ m_{{G}\mc}\circ (\kappa_{{G}\mc} \ot\kappa_{{G}{\mc}}) = m\circ (\rho_C\ot\rho_C)$.
So 
$(\epsilon_\mc \circ m_{{G}\mc}) \circ (\kappa_{{G}\mc} \ot\kappa_{{G}{\mc}}) = (m\circ(\epsilon_\mc\ot\epsilon_\mc))\circ   (\kappa_{{G}\mc} \ot\kappa_{{G}{\mc}}).$

 Since $\kappa_{{G}\mc} \ot\kappa_{{G}{\mc}}$ is an epimorphism by condition 2,
we have $\epsilon_\mc\circ m_{{G}{\mc}}= m\circ (\epsilon_\mc\ot \epsilon_\mc)$; 
that is, $\epsilon_\mc$ preserves the multiplication as required. 
 
 For showing that $\epsilon_\mc$ preserves the unit one needs to show that the lower triangle of Diagram (\ref{eqn:eps}) commutes. Since the outer frame and the {curved} cells do so (see above) and $\rho_I$ is an isomorphism by normality of $F$ and $G$, this follows from commutativity of the central cell, which  is clear by Definition \ref{def:indm}.
 \begin{equation}\label{eqn:eps}
\begin{minipage}{6cm}
\xymatrix@=2.5em{
LGI     \ar[rrr]^{LGe }\ar@/_2pc/@{->}[dd]_{\rho_I} 
&&&LGC \ar[d]_{ \kappa_{G\mc }}\ar@/^2pc/@{->}[dd]^{\rho_C} \\
LI \ar[d]^{ \lambda}  \ar[u]_{L\gamma}    
&& &FG\mc\ar[d]_{\epsilon_\mc}\\
 I\ar[rrru]^{e^\bullet}\ar[rrr]^{e}&&&C
}
\end{minipage}
\end{equation}
\end{Proof}

\begin{Lemma}\label{lem:ad}
For every $\BD$-monoid $\ma$ the $\BD$-morphism $\eta_\ma:=  A\xra{ \lambda_A}GLA\xra{ G\kappa_\ma}G\ma^\circ$  is a monoid morphism $\ma\rightarrow {G}\ma^\bullet$.
\end{Lemma}
\begin{Proof}
Commutativity of diagram (\ref{diag:eta}) shows, that
$\eta_\ma$ preserves multiplication. 
Also, $\eta_A$ respects units as is seen in a similar way.
\end{Proof}

\begin{Lemma}
The families $(\eta_\ma)$ and $(\epsilon_\mc)$ are natural transformations and satisfy the triangle equalities. Thus, $(-)^\bullet$ is  a left adjoint of ${G}$.
\end{Lemma}

\begin{Proof}
Naturality of $\eta$ is clear;   since each $\kappa_{{G}\mc}$ is an epimorphism, $\epsilon$ is natural as well.

Since the underlying functors are faithful it suffices to show that, for any $\BC$-monoid $\mc$ and any $\BD$-monoid $\ma$ respectively the following identities hold
\begin{eqnarray}
{G}\mc \xra{\eta_{{G}\mc} }{G}({G}\mc)^\bullet 
\xra{ {G}\epsilon_{\mc}} {G}\mc &= & id_{GC}\label{1}\\
\ma^\circ \xra{ \eta_\ma^\circ}({G}(\ma^\bullet))^\circ 
\xra{ \epsilon_{\ma^\bullet}}\ma^\circ & = &id_{\ma^\circ}\label{2}
\end{eqnarray}
By definition we have
\begin{eqnarray*}
GC\xra{\eta_{{G}\mc} }{G}({G}\mc)^\bullet &= &GC\xra{ \lambda_{GC}}GLGC\xra{ G\kappa_{{G}\mc}}G(GC)^\bullet\\
GLGC\xra{G\kappa_{{G}\mc} } G({G}\mc)^\circ  \xra{G\epsilon_\mc} GC 
&=& GLGC\xra{ G\rho_C}GC
\end{eqnarray*}
Hence, equation \eqref{1} follows from the first triangle equality for the adjunction $L\vdash G$.
\\
Similarly we have
\begin{eqnarray*}
LA\xra{\kappa_\ma} \ma^\circ\xrightarrow{\eta_\ma^\circ} (G\ma^\bullet)^\circ &=&
LA\xra{ L\lambda_\ma} LGLA\xra{ LG\kappa_\ma}LG\ma^\bullet 
\xra{ \kappa_{G\ma^\bullet}}(G\ma^\bullet)^\circ
\\
LG\ma^\circ\xra{\rho_{\ma^\circ} }\ma^\circ &=&
LG\ma^\bullet\xra{ \kappa_{G\ma^\bullet}}
(G\ma^\bullet)^\circ\xra{  \epsilon_{\ma^\bullet}} \ma^\circ 
\end{eqnarray*}
Since $\rho$ is natural, this implies $\epsilon_{\ma^\bullet}\circ \eta_\ma^\circ\circ\kappa_\ma = \kappa_\ma\circ \rho_{LA}\circ L\lambda_A$.  
Thus, equation \eqref{2} follows by the second triangle equation for $L\dashv G$ since $\kappa_\ma$ is an epimorphism. 
\end{Proof}

\begin{Remarks} 
For the semi-dualization functor $G$:
\begin{enumerate}
\item the basic situation of Definition~\ref{def:BS} is (a) of \cite[Lemma 6.0.1]{Sweedler},
\item Theorem~\ref{prop:leftadjointgeneral} item 2 is trivial in the field case in \cite{Sweedler} and needs specific attention for more general rings (see \cite{Cao}),
\item the minimality condition Theorem~\ref{prop:leftadjointgeneral} item 1(b) is the equivalence of 1) and 3) in \cite[Prop. 6.0.3]{Sweedler}. 
\end{enumerate}
 \end{Remarks}

\subsection{Constructing dual comonoid functors}\label{ssec:dualcomonoid}

Revisiting Sweedler's construction of the left adjoint to the dual algebra functor with a monomorphic canonical transformation $\kappa$, 
we will now analyze the conditions on a symmetric monoidal category 
which make this work. The strategy clearly is to first construct a  basic situation $((-)^\circ,\kappa)$ with respect to the  semi-dualization functor $D\colon\BD^\op\rightarrow\BD$   of a symmetric monoidal closed category\footnote{From now on we write $A^\ast$ instead of $DA$ and $\ma^\bullet$ instead of $\mathsf{Sw}\ma$, as in case of modules.}, that is,  with $(-)^\circ\colon \Mon\BD\lra \BD^\op$, such that Theorem \ref{prop:leftadjointgeneral} can be applied.
In a second step we explain how Sweedler's results and its generalizations 
(see \cite{Cao}, \cite{Wis}) follow. 

\subsubsection{Basic situations}

To gain insight on how, for a $\BD$-monoid $\ma$, the object $\ma^\circ$ 
should be chosen as a subobject of $A^\ast$ in order to obtain a Sweedler functor 
supporting a left adjoint $(-)^\bullet$ of the dual monoid functor, we start with the
 observation that, by Proposition \ref{prop:restr_fct}, 
 there exists a subcategory $\BS$ of $\BD$, closed under the monoidal structure,  
 such that $\mk^\circ = K^\ast$, for each $\BD$-monoid $\mk$ with $K$ in $\BS$.   
For example (see Example \ref{fact:exs} (1)), we could choose $\BS = {_{rd}\BD}$.

Then, for every monoid morphism
 $\ma\xrightarrow{f}\mk_f$ where the codomain $K_f$ of $f$ belongs to $\BS$,  
 the {morphism} $f^\ast$  factors as  $ K_f^\ast = \mk_f^\circ \xrightarrow{f^\ast} A^\ast = K_f^\ast \xrightarrow{{f^\circ}} \ma^\circ\xra{\kappa_\ma} A^\ast$.
The obvious choice for $\ma^\circ$ thus would be the smallest subobject of 
$A^\ast$ through which all such $f^\ast$ factor, and this for $\BS$ as large as possible. See Section \ref{choice} for some remarks concerning the choice of such $\BS$.

Such  factorizations exist, provided that  $\BD$ has
 (extremal episink, mono)-factoriza\-tions in the sense of \cite{AHS}; this is true of $\Mod_R$, for each (commutative) ring $R$.
 
Let $\BS$ be an arbitrary full subcategory of $\BD$.  
Consider, for every monoid $\ma\in\Mon\BD$, the family  of all monoid morphisms $f\colon \ma\rightarrow\mk_f$ where  $\mk_f = (K_f,m_f,u_f)$ has $K_f$ in  $\BS$, and denote by $\mathcal{S}_\ma$ the family of the duals $f^\ast\colon K_f^\ast\rightarrow A^\ast$. 
Form the (extremal episink, mono)-factorization of $\mathcal{S}_\ma$ as
$$ K^\ast_f\xrightarrow{f^\ast} A^\ast = K^\ast_f\xrightarrow{e_f} \ma^\circ\xrightarrow{\kappa_\ma}A^\ast$$
If $h\colon\ma\rightarrow\mb$ is a monoid morphism then the family of all $h^\ast\circ g^\ast$ with $g^\ast\in \mathcal{S}_\mb$ is a subfamily of $\mathcal{S}_\ma$ 
and thus factors as
$ K^\ast_{g\circ h}\xrightarrow{ g^\ast}B^\ast\xrightarrow{h^\ast} A^\ast = K^\ast_{g\circ h}\xrightarrow{e_{g\circ h}} \ma^\circ\xrightarrow{\kappa_\ma}A^\ast.$
With the (extremal episink, mono)-factorization 
$ K^\ast_g\xrightarrow{g^\ast} A^\ast = K^\ast_g\xrightarrow{s_g} \mb^\circ\xrightarrow{\kappa_\mb}B^\ast$
 of $\mathcal{S}_\mb$, we then have the commutative diagram, which has a (unique) diagonal $h^\circ$.
$$\begin{minipage}{6cm}
\xymatrix@=2em{
    K^\ast_g \ar[rr]^{s_g  }\ar[dd]_{ e_{g\circ h} }&&{\mb^\circ } \ar[d]^{\kappa_\mb }\ar@{.>}[ddll]^{h^\circ}  \\
   & &B^\ast\ar[d]^{h^\ast}\\ 
 \ma^\circ  \ar[rr]_{\kappa_\ma}          &&    A^\ast 
}
\end{minipage}$$
Note that, in particular, 
\begin{equation}\label{eqn:fact2}
\mk^\circ\xra{ f^\circ} \ma^\circ  = \mk^\ast \xra{ e_f}\ma^\circ \ .
\end{equation}
This defines a functor $\Mon\BD\xrightarrow{(-)^\circ} \BD^\op$ such that the family  of all $\kappa_\ma$ is a monomorphic natural transformation $(-)^\circ\Rightarrow (-)^\ast\circ |-|$. We, thus, have got a Sweedler functor, depending on $\BS$. Note that one has $\ma^\circ = \ma^\ast$ if and only if $A\in\BS$.\footnote{More precisely one should have said instead of $\ma^\circ = \ma^\ast$ that $\kappa_\ma$ is an isomorphism.}

\begin{Remark}\label{rem:SF}\rm
It is easy to see that this construction of $\ma^\circ$, if applied in the category $\Vect$ to the subcategory $\BS = {_{fd}\Vect}$ gives precisely Sweedler's definition of $\ma^\circ$, i.e., the submodule of $A^\ast$ consisting of all linear forms on $A$, whose kernel contains a cofinite ideal. 
As shown by Sweedler this construction then allows a lift to the desired left adjoint. We will add a remark towards the naturality of his construction in Section~\ref{choice}.

Doing the same construction more generally in $\Mod_R$ with $\BS = {_{fg}\Mod_R}$, for some commutative ring $R$, one certainly obtains a Sweedler functor. There is, however, no reason to assume that this can be lifted in general.
We are going to investigate now, how to find conditions on a commutative ring $R$ and a subcategory $\BS$ of $\Mod_R$ to make this happen. 
\end{Remark}

In view of Remark \ref{rem:Xi}, the first statement of the following lemma must hold, if the basic situation under consideration {is to} lift to a left adjoint.
\begin{Lemma}\label{lem:crux}
Let $\BS$ be a monoidal subcategory of $\BD$.  
Then there exists a natural family $\Phi_{\ma ,\mb} \colon \ma^\circ \ot \mb^\circ \rightarrow (\ma\ot \mb)^\circ$  such that the following diagram commutes.
$$\begin{minipage}{6cm}
\xymatrix@=2.5em{
  \ma^\circ \ot \mb^\circ   \ar[r]^{ \Phi_{\ma ,\mb} }\ar[d]_{\kappa_\ma\ot \kappa_\mb  }&{ (\ma\ot \mb)^\circ} \ar[d]^{\kappa_{\ma\ot\mb} }  \\
 A^\ast\ot B^\ast   \ar[r]_{\Lambda_{A,B}}          &   (A \ot B)^\ast
}
\end{minipage}$$
The morphisms $\Phi_{\ma ,\mb}$ are isomorphisms, provided that $\BS$ is contained in $\BD_D$ and  for each pair  $(\ma,\mb)$ of $\BD$-monoids  the morphism $\Lambda_{A,B}\circ (\kappa_\ma\ot\kappa_\mb)$ is a monomorphism.
\end{Lemma}

\begin{Proof}
Consider the following commutative diagram, where $f$ and $g$ are   all monoid morphisms into $\ma$ and $\mb$ respectively, whose codomains belong to $\BS$.
\begin{equation}\label{eqn:diag}
\begin{aligned}
\xymatrix@=3em{
 K_f^\ast \ot K_g^\ast    \ar[r]^{ e_f\ot e_g }\ar[d]_{\Lambda_{K_f,K_g}  }\ar[dr]^{f^\ast\ot g^\ast}&{ \ma^\circ \ot \mb^\circ} \ar[d]^{ \kappa_\ma\ot \kappa_\mb}  \\
 (K_f\ot K_g)^\ast\ar[d]_{e_{f \ot g}}\ar[dr]^{(f\ot g)^\ast} & A^\ast\ot B^\ast\ar[d]^{\Lambda_{A,B}}\\
 (\ma\ot\mb)^\circ   \ar[r]_{ \kappa_{\ma\ot\mb}}          &   (A\ot B)^\ast
}
\end{aligned}
\end{equation}
The required map would be a diagonal, and this exists, since the family of all $e_f\ot e_g$ is an extremal episink, using our assumption that $\BS$ is closed under the monoidal structure.

For the same reason the set of all $f\ot g$ contains the family of all monoid morphisms $\mk_h\xrightarrow{h} \ma\ot\mb$. Consequently, if $K_h^\ast \xra{ h^\ast}(A\ot B)^\ast = K_h^\ast \xra{ s_h} (\ma\ot\mb)^\circ \xra{ \kappa_{\ma\ot\mb}} (A\ot B)^\ast$ is the (extremal episink, mono)-factorization of the family $h^\ast$, 
the family $(e_{f\ot g})$ contains the family $(s_h)$ and is thus an extremal episink. 
Replacing now in the diagram above the maps $\Lambda_{K_f,K_g}$ by their inverses one obtains again a commutative diagram, which has a diagonal $\iota_{\ma,\mb}\colon (\ma\ot \mb)^\circ \rightarrow \ma^\circ\ot\mb^\circ$ as well. By uniqueness of diagonals this is an inverse of $\Phi_{\ma,\mb}$.
\end{Proof}

There is an alternative description of $\ma^\circ$ in the case of $\BD = \Mod_R$ when $\BS$ is the category of finitely generated modules, provided that $R$ is a noetherian ring. For this we first recall the following fact.

 \begin{Fact}\label{fact:modules}\rm
Let $\ma =(A,m,e)$ be an $R$-algebra.  Then the object $A^\ast$ becomes an $\ma$-bimodule by the operations $A\ot A^\ast \xra{ l}A^\ast$ and $A^\ast\ot A \xra{ r}A^\ast$ characterized by commutativity of the diagrams
\begin{equation}\label{diag:bimod}
\begin{aligned}
\xymatrix@=2.5em{  A\ot A^\ast \ar[r]^{ ev} & I    \\
   A\ot A\ot A^\ast \ar[r]_{ m\ot A^\ast }\ar[u]^{ A\ot l}& A\ot A^\ast\ar[u]_{ev}
}\qquad
\xymatrix@=2.5em{  A^\ast\ot A \ar[r]^{ ev} & I    \\
 A^\ast \ot   A\ot A\ar[r]_{ A^\ast\ot m }\ar[u]^{ r\ot A}& A\ot A^\ast\ar[u]_{ev} \ .}
 \end{aligned}
\end{equation}
As usual we write $af = l(a\ot  f)$ and $fb = r(f\ot b)$.

For every morphism $f\in A^\ast$ there is the algebra morphism $\bar{f}\colon R\rightarrow A^\ast$ with $\bar{f}1 = f$. The image of $A\ot R \xra{ A\ot \bar{f}} A\ot A^\ast \xra{ l}A^\ast$ 
is easily seen to be  a left $\ma$-submodule $\ma f$ of $A^\ast$. We then have
the corresponding algebra morphism $\ma f  \xra{ l_f} \mathrm{End}(\ma_f)$. 
Denoting the elements of $\ma f$ simply by $af$, one sees  that, for $a,b\in A$, $af = bf$  implies $f(a) = f(b)$. Hence, there is a linear form $\zeta$ on $\mathrm{End}(\ma f)$ acting on $\ma f \xra{h}\ma f$  by $\zeta h = f(a)$   if $h(f)  = af$. One now immediately concludes $l_f^\ast(\zeta) = f$.

$f\ma$ is defined analogously.
\end{Fact}

Now the following holds.
\begin{Lemma}\label{lem:Af}
Let $R$ be a  commutative noetherian ring and let $((-)^\circ,\kappa)$ be the basic situation determined by the subcategory $\BS$ of all finitely generated $R$-modules. Then one has for any $R$-algebra $\ma$ and $f\in A^\ast$ 
$$  f\in\ma^\circ \iff  
\ma f \text{ is finitely generated \ }\iff f\ma\text{ is finitely generated.}$$ 
\end{Lemma}
\begin{Proof}
If $\ma f$ is finitely generated then so is  the endomorphism module of $\ma f$, since  $R$ is noetherian\footnote{ $M \simeq R^n/U$ with $n\in \N$ and a submodule	 $U$ in $R$ implies that ${\mathrm{End}}(M)$ is isomorphic to the submodule $\{f \in Hom(R^n,M)\mid U\subset \ker f \}$  of $\hom(R^n,M) \simeq M^n$ and, thus, is finitely generated, since $R$ is noetherian.},  and so then is the image of  the algebra morphism $l_f$. 
Thus, $l_f^\ast$ belongs to the family $\mathcal{S}_\ma$ defining $\ma^\circ$ (see Section \ref{ssec:dualcomonoid}); since $f$ belongs to the image of $l_f^\ast$ by the final statement of Fact \ref{fact:modules}, we get $f\in \ma^\circ$ as required.
The remaining implication  is the implication $(iii)\implies (i)$ of \cite[Prop. 2.6]{Wis}.

The second equivalence follows from the {first} by the  observations that $\ma^\circ =   \ma_{\mathsf{d}}^\circ$ and $f\ma = \ma_{\mathsf{d}}f$, where $\ma_{\mathsf{d}}$ denotes the opposite algebra of $\ma$.
\end{Proof}

\subsubsection{Comultiplications}
In order to obtain a liftable basic situation, we observe first that,  by Proposition \ref{prop:normal}, the following lemma holds.
\begin{Lemma}
Let  $\BS$ be contained in $\BD_D$. Then, for every $\BD$-monoid $\mk =(K,m,u)$ with $K\in \BS$, hence $\mk^\circ = K^\ast$,  the triple $(\mk^\circ, \Lambda_{K,K}^{-1}\circ m^\ast, \lambda^{-1}\circ u^\ast)$ is a comonoid.   
\end{Lemma}

The required comultiplication $m^\bullet$ on $\ma^\circ$  should make diagram \eqref{diag:SW4} commute and $f^\circ$ a monoid morphism. 
This requires (see Equation~\eqref{eqn:fact2} above) that 
$\mu_\ma$ is a diagonal in the diagram below, where $m_f$ and $m_\ma$ denote the respective multiplications. 
It is commutative,  since $f$ is a monoid morphism and $\Lambda$ is natural.
$$\begin{minipage}{6cm}
\xymatrix@=2.5em{
   K_f^\ast  \ar[rr]^{e_f  }\ar[d]_{ m_f^\ast }\ar@{.>}[drr]^{f^\ast}  
   &&{ \ma^\circ} \ar[d]^{\kappa_\ma }  \\
   (K_f\ot K_f)^\ast\ar[d]^{\Lambda_{K_f,K_f}^{-1}}\ar@{.>}[ddrr]^{(f \ot f)^\ast} && A^\ast\ar[dd]^{m_\ma^\ast} \\
   K_f^\ast\ot K_f^\ast\ar[d]_{e_f\ot e_f}\ar@{.>}[dr]^{f^\ast\ot f^\ast} &&  \\
  \ma^\circ\ot \ma^\circ  \ar[r]_{ \kappa_\ma\ot \kappa_\ma}          &    A^\ast\ot A^\ast\ar[r]_{\Lambda_{A,A}} & (A\ot A)^\ast
}
\end{minipage}$$
Thus, $m^\bullet$ exists as required, provided  that $\Lambda_{A,A}\circ (\kappa_\ma\ot \kappa_\ma)$ is a monomorphism. 

Similarly, the diagram below commutes, since $f$ preserves units, and, thus, has a unique diagonal $e^\bullet\colon \ma^\circ \lra I$.
$$\begin{minipage}{6cm}
\xymatrix@=2em{
   K_f^\ast  \ar[rr]^{ e_f }\ar[d]_{u_f^\ast  }\ar[drr]^{f^\ast}& &{\ma^\circ } \ar[d]^{\kappa_\ma }  \\
 I^\ast \ar[d]_{\lambda^{-1}} && A^\ast\ar[d]^{e^\ast}\\
 I   \ar[rr]_{ \lambda}    &      &    I^\ast 
}
\end{minipage}$$
One now gets  by Lemma \ref{lem:Abullet} the following result which   applies, in particular, to any basic situation given by a  full monoidal subcategory $\BS$ of $\BD$ contained in $\BD_D$.
\begin{Lemma}\label{prop:dualcomonoid}
Let $\BD$ be  a  {locally presentable }Êsymmetric monoidal closed category and $((-)^\circ,\kappa)$ a basic situation for its semi-dualization functor, such that each $\ma^\circ$ admits a comultiplication $m$ and a counit $e$ induced from $\ma$ by $\kappa_\ma$. 
Then the basic situation is liftable,  provided that all of the following morphisms are monomorphisms: 
 \begin{enumerate}
\item  $\Lambda_{A,A}\circ (\kappa_{\ma}\ot \kappa_{\ma})$,
\item  $\Lambda_{A,A,A}\circ (\kappa_{\ma}\ot \kappa_{{\ma}}\ot \kappa_{ {\ma}})$.
\end{enumerate}
\end{Lemma}

The essential results of   \cite{Cao} then are that,  given a Dedekind domain $R$ and choosing the class of all finitely generated modules as the category $\BS$ in $\Mod_R$, that $\BS$ consists of $G$-strong modules and that conditions 1 and 2 are satisfied. The stronger statement
  \cite[Lemma 3.b)]{Cao}, that all maps $\Phi_{\ma,\mb}$ are isomorphisms, is not needed, as our arguments show. Indeed, this stronger statement is a consequence as shown in Lemma \ref{lem:crux}.

  The shortcoming of Lemma \ref{prop:dualcomonoid} is the assumption, that all objects of $\BS$ need to be $D$-strong. In \cite{Wis} therefore, for  a noetherian ring $R$,  a different approach for supplying $\ma^\circ$ with a comultiplication   is proposed, a simplified version of  which we will sketch  now for the following reasons: (a) The algebra behind this approach is useful for applying Lemma \ref{prop:dualcomonoid}, and (b) we will show how to improve the main results of \cite{Wis} with respect to functoriality of their construction.

\begin{Lemma}
Let $R$ be a  commutative noetherian ring and let $((-)^\circ,\kappa)$ be the basic situation determined by the subcategory $\BS$ of all finitely generated $R$-modules. Then, for any $R$-algebra $\ma =(A,m,e)$ and every $f\in A^\circ$
$$m^\ast(f) = R^t\circ R^m(f)\in \text{Im}\Pi'\cap \text{\tr{Im}}\Pi'',$$
where  $\Pi ' = R^A\ot \ma^\circ\xra{\id\ot(\theta_\ma\circ\kappa_\ma)}R^A\ot R^A\xra{\Pi_{\ma,\ma}}
  R^{A\times A}$ and 
  $\Pi '' = \ma^\circ\ot R^A\xra{(\theta_\ma\circ\kappa_\ma)\ot \id}R^A\ot R^A\xra{\Pi_{\ma,\ma}}
  R^{A\times A}$.
\end{Lemma}
\begin{Proof} 
Assume $a,b\in A$. Then $R^t\circ R^m(f)(a,b) = f(m(a\ot b)) = bf(a)$.

Since $\ma f$ is finitely generated by Lemma \ref{lem:Af}, there are $g_1,\ldots g_n\in \ma f$, such that  $bf = \sum r_ig_i$. \tr{Choose} $h_1,\ldots h_n\in R^A$ with $h_i(b) = r_i$. Then $bf(a) = \sum f_i(b)g_i(a) = \Pi'(\sum f_i\ot g_i)(a,b)$.  $R^t\circ R^m(f)(a,b)\in \text{im}\Pi''$ follows symmetrically.
\end{Proof}

Now denote the map $\ma^\circ\ot \ma^\circ\xra{(\theta_\ma\circ\kappa_\ma)\ot(\theta_\ma\circ\kappa_\ma)}R^A\ot R^A\xra{\Pi_{\ma,\ma}} R^{A\times A}$
by $\pi$ and assume that
$\theta_\ma\circ\kappa_\ma$ is injective. Then $\text{Im}\Pi'\cap \text{im}\Pi'' = \Pi_{A,A}[(\ma^\circ\ot R^A)\cap (R^A\ot\ma^\circ)] = \pi[\ma^\circ\ot\ma^\circ]$, and one obtains a linear map $\pi[\ma^\circ\ot\ma^\circ]\xra{\phi}\ma^\circ\ot\ma^\circ$ by taking preimages, since $\pi$ is injective. With $m^\bullet = \ma^\circ\xra{}\pi[\ma^\circ\ot\ma^\circ]\xra{\phi}\ma^\circ\ot\ma^\circ$ we get (see also \cite{Wis})
\begin{Corollary}
\label{lem:wisind}
Let $R$ be a noetherian ring and $\ma = (A,m,e)$ an $R$-algebra.  
Then there exists an $R$-linear map $\ma^\circ \xra{m^\bullet}\ma^\circ\ot \ma^\circ$ such that the outer frame of Diagram (\ref{diag:wis1}) commutes (that is by  Section \ref{ssec:struct}, that $m^\bullet$ is induced from $\ma$ by $\kappa_\ma$), provided that $\ma^\circ$ is pure in $R^A$. 
\end{Corollary}

Obviously, if $\ma^\circ$ is a pure submodule of $R^A$ with embedding  $\theta_\ma\circ\kappa_\ma$, then $\ma^\circ$ is a pure submodule of $A^\ast$ via $\kappa_\ma$ and, consequently, the morphism $\kappa_\ma\ot\kappa_\ma$ is a  monomorphism.
Commutativity of the diagram below, with $\Phi$ as in Lemma \ref{lem:crux}, now shows  
 that then $\Lambda_{A,A}\circ (\kappa_\ma\ot\kappa_\ma)$ is a monomorphism as well, provided that $\Pi_{A,A}$ is a monomorphism. Similarly, the morphism $\Lambda_{A,A,A}\circ (\kappa_\ma\ot\kappa_\ma\ot\kappa_\ma)$ is a monomorphism, provided that $\Pi_{A,A,A}$ is monic. 
\begin{equation}\label{eqn:wis}
 \begin{aligned}
\xymatrix@=3em{
\ma^\circ\ot  \ma^\circ \ar[r]^{\kappa_\ma\ot\kappa_\ma} \ar[d]^{\Phi_{\ma,\ma}}  & A^\ast\ot A^\ast   \ar[rr]^{\theta_\ma\ot\theta_\ma}\ar[d]_{\Lambda_{A,A}} &&  R^A\ot R^A \ar[d]^{\Pi_{\ma,\ma}}\\
(\ma\ot \ma)^\circ \ar[r]^{\kappa_{A \ot A}} & (A\ot A)^\ast   \ar[r]^{\theta_{A\ot A}}& R^{A\ot A}\ar[r]^{R^{t_A}}& R^{A\times A} 
}
    \end{aligned}
\end{equation}

We thus may, by Lemma \ref{prop:dualcomonoid} and the simple fact, that all maps $\Pi_{A,A}$ and $\Pi_{A,A,A}$ are injective, if $R$ is noetherian (see \cite{Wis}), strengthen the lemma above as follows.
\begin{Proposition}\label{prop:pure}
Let $R$ be a noetherian ring such that, for every $R$-algebra $\ma$, the module $\ma^\circ$ is a pure submodule of $R^A$.
Then the basic situation determined by the class of finitely generated $R$-modules is liftable to a dual coalgebra functor $(-)^\bullet\colon\Alg_R\lra\Coalg_R^\op$.
\end{Proposition}
 
\subsection{The minimality condition}
 
  The remaining condition of Theorem \ref{prop:leftadjointgeneral} to be satisfied now is, that 
 $\ma^\bullet$ not only be a subcomonoid of $A^\ast$ induced by $\ma$, but the smallest such.

 We will discuss this only in the case where $\BD$ is  $\Mod_R$,  the category of modules over a commutative ring $R$. 
 Assuming here again that the maps $\psi:=\Lambda_{A,A}\circ (\kappa_\ma\ot\kappa_\ma)$ are injections, the minimality requirement amounts to saying that $\ma^\circ$ is the preimage under $m^\ast$ of $X:=Im\psi$.
This is known to be true over a field (see \cite[Prop. 6.0.3 ]{Swe}).

We review  the core of Sweedler's argument as follows, where we first
 recall that a comonoid $(C, {\mu}, {\epsilon})$ is a
 subcomonoid of $A^\ast$ induced from $\ma$ by $\kappa_\ma$, provided that the following diagram	commutes,
where  $\mathsf{Sw}(\ma) = (C, {\mu}, {\epsilon})$. 
\begin{equation}\label{diag:pullback}
\begin{aligned}
\xymatrix@=2em{
    C \ar[d]_{ \kappa_\ma }\ar[r]^{  {\mu}}& C\ot C\ar[r]^{\kappa_\ma\ot \kappa_\ma} & A^\ast\ot A^\ast  \ar[d]^{\Lambda_{A,A } }     \\
A^\ast\ar[rr]^{m^\ast}
   &  &    (A\ot A)^\ast 
}
\end{aligned}
\end{equation}
 
\begin{Lemma}\label{lem:Sweedler}
Let $R$ be a commutative noetherian ring. 
Let $\BS$ be a monoidal subcategory of $\Mod_R$, such that the basic situation determined by $\BS$ is liftable.
Then every subcoalgebra $\mc$ of $A^\ast$ induced by $\ma$ is a subcoalgebra of $\ma^\bullet$, provided that 
$\BS$ contains all finitely generated modules.
\end{Lemma}

\begin{Proof} 
By assumption	we may assume  $C\subset A^\ast$ and that
 diagram~\eqref{diag:pullback} commutes.
 Thus, with notation as above, we have
 $m^\ast f\in X = Im\psi$ for each $f\in C$. It thus suffices to prove the implications  $m^\ast f\in X \implies\ma f\in {_{fg}\Mod_R} \implies   f\in \ma^\circ$.

Now $m^\ast f \in X$ means $\psi(m^\ast f) = \sum_ig_i\ot h_i$ with $g_i,h_i\in\ma^\circ\ot\ma^\circ$. 
With $a,b\in A$, we have $af(b) = f(m(b\ot a)) =  m^\ast f(b \ot a)$, hence $\psi(m^\ast f)(a\ot b) =  (\sum_ig_i\ot h_i)(b\ot a) = \sum_ig_ib\cdot h_ia = (\sum_i(h_ia)g_i)(b).$
The map $af$ (that is, an arbitrary element of $\ma f$) is thus contained in the subspace 
of $A^\ast$ generated by the maps $g_i$. 
Thus, $\ma f$ is finitely generated, and this proves the claim by  Lemma \ref{lem:Af}.
 \end{Proof}
 
Thus, $\ma^\bullet$ is the largest subcoalgebra of $A^\ast$ induced by $\ma$, and by Theorem~\ref{prop:leftadjointgeneral} we deduce:

\begin{Proposition}\label{fact:Sweedler}
Let $R$ be a noetherian ring. If  the basic situation determined by $\BS = {_{fg}\Mod_R}$ is liftable to a dual coalgebra functor $(-)^\bullet\colon\Alg_R\rightarrow \Coalg_R^\op$, 
 then this functor is left adjoint to the dual algebra functor and hence is the generalized Sweedler dual $\mathsf{Sw}$.
 \end{Proposition}

\section{The  generalized Sweedler dual of an $R$-algebra}\label{sec:fin}

By the results of the last section it remains 
 to find conditions on a noetherian commutative ring $R$ such that  the maps $\Lambda_{A,A}\circ (\kappa_\ma\ot\kappa_\ma) $ and $\Lambda_{A,A,A}\circ (\kappa_\ma\ot\kappa_\ma\ot\kappa_\ma)$ are injective for each $R$-algebra $\ma$. And this will be the case (except for the trivial case of an absolutely flat noetherian ring), if for every $R$-algebra $\ma$ the module $\ma^\circ$, defined by the  basic situation determined by $\BS = {_{fg}\Mod_R}$, is a pure submodule of $R^A$.

Thus, the following algebraic observation, whose proof is quite simple, is useful.
\begin{Lemma}[\cite{Wis}]\label{lem:Wis} 
Let $R$ be a hereditary noetherian ring $R$. Then, for  every $R$-algebra $\ma$, 
 $A^\circ$ is a pure submodule of $R^A$. 
\end{Lemma}

Now the following extensions of the main results of \cite{Wis} and \cite{Cao} are easy consequences.
\begin{Theorem}\label{corr:Ded}
Let $R$ be a noetherian ring. Then the following hold.
\begin{enumerate}
\item If $\ma$ is an $R$-algebra, where $\ma^\circ$ is a pure submodule of $R^A$, then
\begin{enumerate}
\item $\ma^\bullet$ is an $R$-coalgebra.
\item If $\ma$ is the underlying algebra of a Hopf algebra then $\ma^\bullet$ is the underlying coalgebra of a Hopf algebra.
\end{enumerate}
\item If $R$, in addition, is hereditary, then the construction of  $\ma^\bullet$ defines a functor $\Alg_R\rightarrow\Coalg_R^\op$ and the following hold. 
\begin{enumerate}
\item This functor is left adjoint to the dual algebra functor and, hence, the generalized Sweedler dual $\mathsf{Sw}$.
\item As the opmonoidal left adjoint of the dual algebra functor 
$\mathsf{Sw}$ is  opmonoidal and so induces a functor $\Bialg_R\lra\Bialg_R^\op$ and, by restriction, a functor 
$\Hopf_R\rightarrow\Hopf_R^\op$.
\item $\mathsf{Sw}$, considered as a functor $\Bialg_R\lra\Bialg_R^\op$, has a right adjoint and thus yields a dual adjunction on $\Bialg_R$. 
\end{enumerate}
\item If $R$ is a Dedekind domain, $\mathsf{Sw}$ is strong opmonoidal.
\item If $R$ is  absolutely flat,  then the dual adjunction of 2(c) can be  restricted to a dual adjunction  on $\Hopf_R$.
\end{enumerate}
\end{Theorem}
\begin{Proof}
Only 1(b) and 4 still need an argument.
By Lemma \ref{lem:crux}, $\Lambda_{\ma,\ma}$ is an isomorphism since $\kappa_\ma\ot\kappa_\ma$ is injective. It now follows dually to the proof of statement 1 of Proposition \ref{prop:normal} that, if $\ma$ is the underlying algebra of a bialgebra, so is $\ma^\bullet$. The same holds for Hopf algebras by the proof of Proposition \ref{prop:hopf}, since $\kappa_\ma$ is injective. Statement 4 is a corollary to Proposition \ref{prop:dualadj2} since the category $\Hopf_R$ of Hopf algebras over $R$ is reflective in $\Bialg_R$, provided that the ring $R$ is absolutely flat (see \cite{HEP_JA}).
\end{Proof}

\section*{Concluding remarks}
\subsection*{The choice of $\BS$}\label{choice}
 As already mentioned at the beginning of Section \ref{ssec:dualcomonoid}, the obvious choice for the definition of $(-)^\circ$ is our factorization method with a suitable $\BS$.
  
 Though, for each ring $R$ and for any $\BS$, one will obtain a Sweedler functor this way, this will not be liftable to a left adjoint of the dual algebra functor in general. {In fact, in order to lift, one needs $\BS$ to contain ${_{f\!gp}\Mod_R}$ by Proposition \ref{prop:restr_fct} in connection with Example \ref{fact:exs}.2. However, the bigger $\BS$ is than this, the bigger than absolutely necessary must be the part of $\Mod_R$ on which the left adjoint would have to be build over $(-)^\ast$, noting that $\ma^\circ = A^\ast \iff A\in\BS$,
which clearly reduces the chances of the construction being possible.}  On the other hand, since our only minimality criterion is Lemma \ref{lem:Sweedler}, one should have $\BS$ containing ${_{f\!g}\Mod_R}$! It thus is well motivated to choose $\BS = {_{f\!g}\Mod_R}$ for \em any \em ring and see how far one gets (note that by the above we don't need all finitely generated modules to be projective - it suffices that they are $D$-strong); and, if every finitely generated $R$-module would be projective (as for a field), this is the obvious choice! 
Thus, Sweedler's definition of $\ma^\circ$ is not really surprising. And that his construction works in $\Vect$, is a special instance of the considerations above, since in this category the additional assumptions 1. and 2. of Lemma \ref{prop:dualcomonoid} are satisfied trivially.
\subsection*{Generalizations}
Most of the results on categories of $R$-modules for a commutative ring $R$ can be generalized to arbitrary finitary varieties with a commutative theory, that is in the language of universal algebra, to entropic varieties. We have refrained  from doing so because we didn't see a way of applying our criteria in categories more general than $\Mod_R$. The paper \cite{Abu} tries to generalize results from \cite{Wis} to categories of semimodules over semi\-rings, and this would have been an example of that kind of generalization; unfortunately, however, there seems to be a gap in its attempt to prove the crucial result, that the maps labeled $\Pi$ in subsection \ref{ssec:struct} are monic for noetherian semirings (see the proof of Proposition \ref{prop:pure} for the importance of this property). The proof of Lemma 5.9 of this paper, on which this would be based, assumes that each semimodule is a directed colimit not only of finitely generated subsemimodules, but of so-called \em uniformly finitely generated \em ones. And we were not able to close this gap. 

\subsection*{Acknowledgment} 
We are grateful to the referee's  thorough work which led to our reorganization of the paper. Crucially, for detecting a false lemma. We are grateful also to Exequiel Rivas for pointing out that Definition 16 of the published version was too weak to prove Proposition 17.

\end{document}